\documentclass[12pt]{amsart}

\usepackage{amssymb}

\usepackage{enumitem}

\usepackage{graphicx}
\usepackage{overpic}
\usepackage{wrapfig}
\usepackage[font=small,labelfont=rm]{caption}

\makeatletter
\@namedef{subjclassname@2020}{%
  \textup{2020} Mathematics Subject Classification}
\makeatother

\newcommand{\bi}[1]{{\bf\itshape #1}}
\newcommand{\R}{{\mathbb R}}		
\newcommand{\C}{{\mathcal C}}           
\newcommand{\cn}{\colon}                

\newcommand{\area}{\mathop{\rm Area}}
\newcommand{\dist}{\mathop{\rm dist}}
\newcommand{\lat}{{\mathcal L}}
\newcommand{\Z}{{\mathbb Z}}
\newcommand{\len}{\mathop{\rm Length}}
\newcommand{\bmat}{\left[\begin{matrix}}
\newcommand{\emat}{\end{matrix}\right]}

\newcommand{\tc}{\operatorname{\pmb\tau}}   

\newcommand{\aff}{\operatorname{Aff}\nolimits}

\newtheorem*{claim}{Claim}
\newtheorem{theorem}{Theorem}[section]
\newtheorem{corollary}[theorem]{Corollary}
\newtheorem{lemma}[theorem]{Lemma}
\newtheorem{proposition}[theorem]{Proposition}

\theoremstyle{definition}
\newtheorem{definition}[theorem]{Definition}
\newtheorem{remark}[theorem]{Remark}

\numberwithin{equation}{section}

\catcode`@=11 \@mparswitchfalse  

\newcounter{mnotecount}  
\renewcommand{\themnotecount}{\arabic{mnotecount}} 
\newcommand{\mnote}[1]
{\protect{\stepcounter{mnotecount}}$^{\mbox{\footnotesize  $
      \bullet$\themnotecount}}$\marginpar{\parbox[b]{1.2in}{\raggedright\tiny\em
	 \themnotecount:\! #1}} }

\begin{document}


\baselineskip=17pt


\title[Bounding the number of integer points near a convex curve]{Bounding the number of
lattice points near a convex curve by curvature}

\author[R. Howard]{Ralph Howard}
\address{Department of Mathematics\\ University of South Carolina\\
Columbia, SC 29208}
\email{howard@math.sc.edu}

\author[O. Trifonov]{Ognian Trifonov}
\address{Department of Mathematics\\ University of South Carolina\\
Columbia, SC 29208}
\email{trifonov@math.sc.edu}

\date{}

\begin{abstract}
We prove explicit
bounds on the number of lattice points on or near a convex curve in terms
of geometric invariants such as length, curvature, and affine arclength. In several of our results we obtain the best possible constants.
Our estimates hold for lattices more general than the usual lattice of
integral points in the plane.
\end{abstract}



\maketitle

\section{Introduction}
Our goal in this paper is to give explicit and as sharp as possible
bounds on the number of lattice points on or near a convex curve in terms
of geometric invariants such as length, curvature, and affine arclength for  lattices more general than the usual lattice of
integral points in the plane.

\begin{definition}
Let $v_0,v_1,v_2\in \R^2$ be vectors with $v_1$ and $v_2$
linearly independent.  Then, the \bi{lattice generated by $v_1$
and $v_2$ with origin $v_0$} is
$$
\lat=\lat(v_0,v_1,v_2) = \{ v_0+m v_1 + nv_2: m,n\in \Z\}.
$$
 \end{definition}

Note that the elements of such a lattice need not have
integral or even rational components.
An invariant of a  lattice is the area spanned by
$v_1$ and $v_2$
$$
A_\lat := | v_1\wedge v_2|
$$
where $v_1\wedge v_2$ is the determinant of the  of the $2 \times 2$
matrix with columns $v_1$ and $v_2$.

If $\C$ is a curve of differentiability class $C^2$ and whose
curvature $\kappa$ is positive,  then the \bi{total curvature}
of $\C$ is
$$
\tc(\C):=\int_\C \kappa\,ds
$$
 where $s$ is arclength along
 $\C$ and the \bi{radius of curvature} of $\C$ is $\rho = 1/\kappa$.
The following are representative of our results.

\begin{theorem}\label{intro1}
Let $\C$ be a $C^2$ curve with total curvature at most $\pi$ and whose radius
of curvature has a lower bound $\rho\ge R$ for some positive constant $R$.
Let $\lat$ be a lattice with
$$
\len(\C) \le 2(A_\lat R)^{1/3}
.$$
Then,   $\C$ contains at most two points of $\lat$.
\end{theorem}

This generalizes a theorem of Schinzel (whose proof first appeared in the paper
\cite[Lemma 2]{Zyg} of Zygmund) where $\C$ is an arc of a circle
  and the
lattice is $\Z^2$.  In 
\cite{Cilleruelo} Cilleruelo
shows that when  $\C$  is an arc of a circle centered at the origin
the sharp form of this inequality has the constant $2$ replaced
by $2 \sqrt[3]{2}$.  In our result, with more general lattices and more general
curves, the constant $2$ is the best possible (see Remark \ref{2_is_best}
below.)

\begin{theorem}\label{intro2}
Let $\C$ be a $C^2$ curve with total curvature $\tc(\C)=\int_\C \kappa\,ds \le \pi$
and whose radius of curvature satisfies $\rho\ge R_1$ for some $R_1>0$.  Then for
any lattice $\lat$
$$
\#(\C\cap \lat) < 2 + \frac{\len(\C)}{(A_\lat R_1)^{1/3}}
$$
If also $\rho\le R_2$, then
$$
\#(\C\cap \lat)\le 2 +  \left( \frac{R_2\tc(\C)}{A_\lat R_1} \right)^{1/3}
\len(\C)^{2/3}
$$
\end{theorem}

This result is close to optimal:

\begin{theorem}\label{intro:example}
Let $\lat$ be a lattice and $n\ge 2$ an integer.  There
is a convex curve $\C$ that contains exactly $n$ points of $\lat$,
and lower and upper bounds
$$
R_1 = \min_{P\in \C} \rho(P), \qquad R_2=\max_{P\in \C} \rho(P)
$$
for the radius of curvature of $\C$, so that both the inequalities
\begin{equation}\label{ineqs:example}
\frac{\len(\C)}{(R_1A_\lat)^{1/3}} < n+2, \qquad \left( \frac{R_2\tc(\C)}{A_\lat R_1} \right)^{1/3}
\len(\C)^{2/3}< n+2
\end{equation}
hold.
\end{theorem}

The foundational result in this subject is the 1926 paper,
\cite{Jar}, of Jarn\'{\i}k who proved that
  the number of integer
points on a strictly convex closed curve of length $L >3$ does not exceed
$3(2\pi)^{-1/3}L^{2/3} + O \left ( L^{1/3} \right )$ and the exponent and the
constant of the leading term are best possible.   Therefore,
the exponent $2/3$ in Theorem \ref{intro2}  is as good as can
be expected.

Using that the affine image of a lattice is a lattice, that every
ellipse is the affine image of a circle, and that affine
arclength (defined in Section \ref{sec:ellipses})
is also invariant under affine maps
we can transfer results about circles to results about ellipses.
One such result is

\begin{theorem}\label{intro2.5}
Let $\C$ be an arc on an ellipse with affine arclength $\aff(\C)$.
Then for any lattice $\lat$
$$
\#(\C\cap \lat) \le 2 + \frac{\aff(\C)}{A_\lat^{1/3}}.
$$
\end{theorem}

We can also estimate the number of points close to a lattice.  This
involves another invariant of a lattice $\lat$, the minimum distance
between any two of its points
$$
d_\lat = \min\{ \| P-Q\| : P,Q\in \lat \text{ and } P \ne Q\}.
$$

\begin{theorem}\label{intro3}
Let $\C$ be a convex arc with total curvature at most $\pi$ with radius of curvature
bounded
 by $R_1\le \rho\le R_2$.  Let
$\lat$ be a lattice and $\delta>0$ with
$$
\delta < \min\left\{ R_1, \ \frac{d_\lat^2}{2( R_2 + d_\lat + \sqrt{(R_2+d_\lat)^2 - d_\lat^2})}
\right\}
$$
and
$$
\frac{A_\lat}{2} - L\delta -\frac32 \delta^2>0.
$$
Then,
$$
\#\{Q \in \lat: \dist(\C,Q)< \delta\}
< 2 + \frac{L}{\big(R_1( A_\lat -2 L\delta - 3 \delta^2)\big)^{1/3} }
$$
where $L = \len(\C)$.
\end{theorem}

Theorems estimating the number of lattice points {\it close} to a curve
are more recent. In 1974 Swinnerton-Dyer improved the exponent in Jarn\'{\i}k's result for curves which are  dilations of a fixed convex $C^3$ curve. In 1989 Huxley \cite{Hux} obtained upper bounds for the number
of lattice points close to the curve $y = f(x)$, $x \in [M,2M]$ assuming $f$
satisfies certain smoothness conditions. In particular, Huxley generalized
Swinnerton-Dyer's result. A number of papers containing new upper bounds for
the number of lattice points close to a curve and applications to different
arithmetic functions ensued. For survey of such estimates and their
applications see the papers  \cite{FT} and \cite{FGT}.

Most of our results are based on some new results on the differential
geometry of plane convex curves which are of interest on their own
right.

\begin{theorem}\label{intro4}
Let $\C$ be a $C^2$ curve with positive curvature and total curvature
$\int_\C \kappa\,ds \le \pi$.  If $\C$ intersects a circle
of radius $R$ in at least $3$ points,  then there is
a point on $\C$ with $\kappa = 1/R$.
\end{theorem}

The structure of this paper is as follows.

Section \ref{sec:lattices} gives  basic facts about lattices
and affine maps.

Section \ref{sec:basic_est} contains basic estimates
we will be using.  The proofs here own a lot to the ideas
in the paper \cite{Cilleruelo-Granville} of Cilleruelo and Granville.

Section \ref{sec:diff_geo} has the proofs of the differential
geometric results we require.

Section \ref{sec:ellipses} starts with results about the number
of points on a circular arc that are on a general lattice $\lat$.
Then the affine invariance of the collection of lattices and affine
arclength under affine maps is used to transfer these results to  the case
of lattice points on an arc of an ellipse.  The results
are new even in the case of the lattice $\lat=\Z^2$.

Section \ref{sec:on_curve} has estimates on the number of
points of a lattice $\lat$ on a convex curve in terms of $A_\lat$
and bounds on the length and curvature of the curve.

Section \ref{sec:near_curve} contains estimates on the number
of points of a lattice $\lat$ within $\delta$ of a convex arc
in terms of $A_\lat$, $d_\lat$, and bounds on the length and
curvature of the curve.

In Section \ref{sec:examples} we show that two of our results are close to being sharp.

\section{Lattices and affine maps.}
\label{sec:lattices}

\begin{definition}
	An \bi{affine map} $\phi\cn \R^n\to\R^n$ is a map
	of the form
	$$
	\phi(v) = Mv+b
	$$
	where $M$ is a non-singular linear map.  Define
	$$
	\det(\phi) = \det(M).
	$$
\end{definition}

The set of lattices is invariant under affine maps.

\begin{proposition}\label{lat-affine-inv}
Let $\phi \cn \R^2 \to \R^2$ be the affine map
$$
\phi(v) = Mv+b.
$$
  Then,
the image of the lattice $\lat(v_0,v_1,v_2)$ under $\phi$ is
$$
\phi\big[\lat(v_0,v_1,v_2)\big] = \lat(\phi(v_0),Mv_1,Mv_2))
$$
and if $\lat=\lat(v_0,v_1,v_2)$ and $\lat^*= \phi\big[\lat]$ is its
image
then
$$
A_{\lat^*} = |\det(\phi)| A_{\lat}.
$$
\end{proposition}

This is straightforward and the proof is left to the reader.

\begin{proposition}\label{A_lat-est}
Let $P_0$, $P_1$, and $P_2$ be three non-collinear points of $\lat
=\lat(v_0,v_1,v_2)$.
Then, the area of the triangle $\triangle P_0P_1P_2$
is an integral multiple of $A_\lat/2$ and therefore
$$
\area( \triangle P_0P_1P_2) \ge \frac 12 A_\lat.
$$
\end{proposition}

\begin{proof}
By the definition of the lattice $\lat$ there are integers $m_j,n_j$
with $0 \le j \le 2$ so that
\begin{align*}
P_j = v_0 + m_jv_1 + n_j v_2.
\end{align*}
Since translation does not change areas,  we can assume $P_0=v_0$.
Then, the area of $\triangle P_0P_1P_2$ is
\begin{align*}
\area( \triangle P_0P_1P_2) &= \frac{1}{2} | (P_1-P_0)\wedge (P_2-P_0)|\\
	&=\frac12 |(m_1v_1 + n_1v_2)\wedge (m_2v_1+ n_2v_2)|\\
	&= |m_1n_2 - m_2n_1| \frac{A_\lat}{2} \\
	&\ge \frac{A_\lat}{2}
\end{align*}
as $|m_1n_2 - m_2n_1|\ge 1$ because it is an integer.
\end{proof}

\section{Conventions and basic geometric estimates.}
\label{sec:basic_est}

All our curves will be of the differentiable class $C^2$
with nonvanishing first and second derivative vectors.  If the
orientation (direction of increasing parameter) of a curve is
is reversed, it changes the sign of the curvature.  As curves
with nonvanishing second derivative have nonvanishing curvature,
by possibly changing the orientation of the curve, we can,
and do, assume all our curves have positive curvature.
If we have a finite set of points $\mathcal F$ on $\C$,
for example if $\# \mathcal F=n$, then we order the points
$\mathcal F = \{ P_1, P_2 ,\ldots, P_n\}$
in the order given by the orientation of the curve.
This implies
that $P_{j+1}$ is between $P_j$ and $P_{j+2}$.

\begin{proposition}\label{prop:Area_pts_circle}
If $\triangle P_0P_1P_2$ is a triangle and its vertices $P_0$, $P_1$, and  $P_2$
are on a circle $\C$ of radius $R$, then the area   of the triangle is
$$
\area(\triangle P_0P_1P_2) = \frac{abc}{4R}
$$
where $a$, $b$ and $c$ are the side lengths of the triangle.  Also, the
area satisfies
the inequality
$$
\area (\triangle P_0P_1P_2) < \frac{(a+b)^3}{16R} .
$$
\end{proposition}

\begin{proof}
The formula for the area is a result attributed to Heron of Alexandria
\cite[Eq.~1.54 p.~13]{Coxeter-Geometry}.
To prove the inequality,  note $c<a+b$ as $a$, $b$ and $c$ are the side lengths
of a triangle.  By the Arithmetic-Geometric mean inequality,
$ab\le (a+b)^2/4$ and therefore
$$
\area (\triangle P_0P_1P_2)=  \frac{abc}{4R} < \frac{\left ( (a+b)^2/4\right )(a+b)}{4R}
= \frac{(a+b)^3}{16R}.
$$
\end{proof}

The next two results are generalizations of results of Cilleruelo and Granville
\cite{Cilleruelo-Granville} from circular arcs to more general curves.  The
proofs are basically axiomatizations of their arguments.

\begin{theorem}[Basic estimate for closed curves]\label{basic_est_closed}
Let $\C$ be a closed curve and $P_1, P_2 ,\ldots, P_N$
points on $\mathcal C$ listed in cyclic order around $\mathcal C$
with the convention $P_{N+1}=P_1$ and $P_{N+2} = P_2$.
Assume there are positive constants $A_0$ and $R_0$ such that
\begin{enumerate}
\item[(a)]  For all $j$
$$
	\frac{A_0}{2} \le  \area(\triangle P_j P_{j+1} P_{j+2})
$$
\item[(b)] For each $j\in \{1,2 ,\ldots,  N\}$ the points $P_j$, $P_{j+1}$, and $P_{j+2}$
	are on a circle of radius $\ge R_0$.
\end{enumerate}
Then
$$
N < \frac{\len(\C)}{(A_0R_0)^{1/3} }.
$$
\end{theorem}

\begin{theorem}[Basic estimate for open curves]\label{basic_est_open}
Let $\C$ be an immersed  curve and $P_1, P_2 ,\ldots, P_N$
points on $\mathcal C$ listed in order along $\mathcal C$.
Assume there are positive constants $A_0$ and $R_0$ such that
\begin{enumerate}
\item[(a)]  For all $j$
	$$
		\frac{A_0}{2} \le  \area(\triangle P_j P_{j+1} P_{j+2})
	$$
\item[(b)]  For each $j$ with $1\le j \le N-2$ the points $P_j$, $P_{j+1}$, and $P_{j+2}$
	are on a circle of radius $\ge R_0$.
\end{enumerate}
Then
$$
N <2+\frac{\len(\C)}{(A_0R_0)^{1/3} }.
$$
\end{theorem}

\begin{proof}[Proof of Theorem \ref{basic_est_closed}]
Let $R_j$ be the radius of the circle through $P_j$, $P_{j+1}$, and $P_{j+2}$.
To simplify notation,  we set
$$
a_j := \| P_{j+1} - P_j\|.
$$
Then by Proposition \ref{prop:Area_pts_circle} and using $R_j\ge R_0$
\begin{align*}
\area(\triangle P_j P_{j+1} P_{j+2})&<
\frac{( a_j+ a_{j+1})^3}{16R_j} \\
&\le \frac{\big( a_j + a_{j+1}\big)^3}{16R_0}
\end{align*}
Combining this with assumption (a) gives
$$
1 < \frac{(a_j+a_{j+1})^3}{8A_0R_0} .
$$
Take cube roots
$$
1<  \frac{ a_j + a_{j+1}}{2(A_0R_0)^{1/3}}
$$
and sum on  $j$
\begin{align*}
	N < \sum_{j=1}^N \frac{a_j + a_{j+1}}{2(A_0R_0)^{1/3}}
	= \frac{1}{(A_0R_0)^{1/3}} \sum_{j=1}^N a_j,
\end{align*}
where we have used $\sum_{j=1}^N a_{j+1}=\sum_{j=1}^N a_j$.
This sum is the length of a polygon inscribed in $\C$ and thus
$\sum_{j=1}^N a_j \le \len(\C)$ which completes the proof.
\end{proof}

\begin{proof}[Proof of Theorem \ref{basic_est_open}]
	As in the proof of Theorem \ref{basic_est_closed} we have
	$$
	1<  \frac{ a_j+ a_{j+1}}{2(A_0R_0)^{1/3}}
	$$
but this time only holding for $1\le j \le N-2$.  Sum on this to
get
\begin{align*}
	N-2&< \sum_{j=1}^{N-2} \frac{a_j+a_{j+1}}{2(A_0R_0)^{1/3}}
= \frac{1}{2(A_0R_0)^{1/3}} \left ( \sum_{j=1}^{N-2}  a_j
+  \sum_{j=1}^{N-2}  a_{j+1}\right )\\
&< \frac{1}{(A_0R_0)^{1/3}} \sum_{j=1}^{N-1}  a_j
\le \frac{\len(\C)}{(A_0R_0)^{1/3}} .
\end{align*}
\end{proof}

\section{Some differential geometry.}
\label{sec:diff_geo}

Let $\C$ be a $C^2$ plane curve and let $\gamma\cn [a,b]\to \R^2$ be a unit speed,
that is $\|\gamma'(s)\|=1$ for all $s$, parametrization of $\C$.  Let
$\mathbf t(s) = \gamma'(s)$ be the unit tangent and $\mathbf n(s)$ the
unit normal where we choose $\mathbf n$ to be $\mathbf t$ rotated by
$\pi/2$ in the positive direction.  Then the \bi{curvature} function along
$\C$ is defined by
$$
\frac{d\mathbf t}{ds} = \kappa(s) \mathbf n.
$$
As remarked above, we orient all our curves so that the curvature is positive.

There is another way to define curvature which will be useful to us.  As
$\mathbf t(s)$ is a unit vector,   it can be written as
$$
\mathbf t = (\cos(\theta(s)), \sin(\theta(s)))
$$
where $\theta$ is a $C^1$ function and is the angle the tangent makes with the
positive $x$-axis.  Then
$$
\frac{d\mathbf t}{ds} = \frac{d\theta}{ds} (-\sin(\theta(s)),\cos(\theta(s)))
= \kappa(s) \mathbf n(s).
$$
Therefore,  the curvature is the rate of change of the angle with respect to arclength:
$$
\kappa = \frac{d\theta}{ds}.
$$
The \bi{total curvature} of $\C$ is the integral of curvature with respect
to arclength and is the total change in the angle of the
tangent vector:
$$
\tc(\C):=\int_\C \kappa\,ds = \int_a^b  \frac{d\theta}{ds} \,ds = \theta(b)-\theta(a).
$$
This interpretation makes it easy to compare the total curvature of two
curves with the same endpoints.

\begin{figure}[h]
\begin{overpic}[height=.75in]{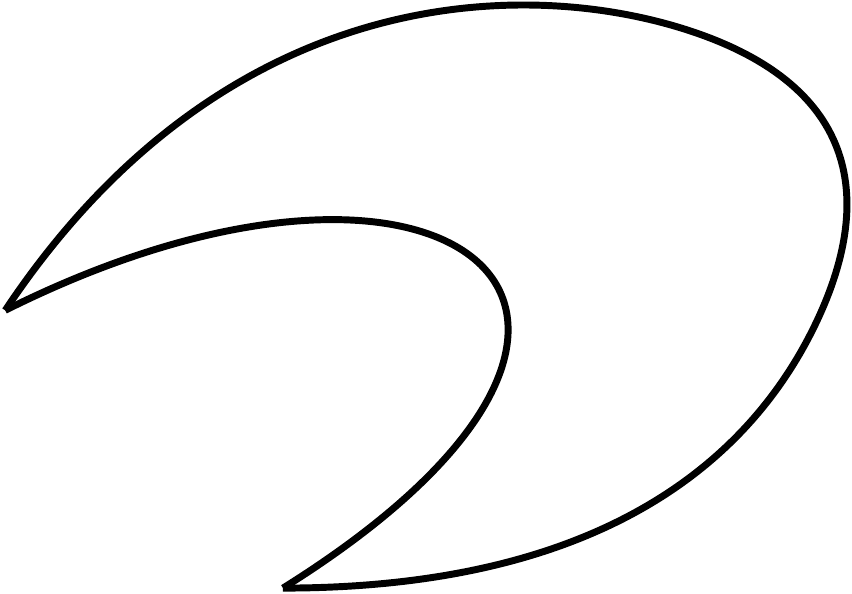}
\put(43,29){$\C_1$}
\put(82,45){$\C_2$}
\end{overpic}
\caption{The total curvature of $\C_2$ is greater than the total
curvature of $\C_1$.}
\label{fig:inside}
\end{figure}

\begin{proposition}\label{prob:comp-total}
Let $\C_1$ and $\C_2$ be convex curves with the same endpoints and with
$\C_1$ inside $\C_2$ in the sense that $\C_1$ is inside the convex hull
of $\C_2$ (see Figure \ref{fig:inside}).  Then  the total curvature of
$\C_1$ is less than or equal to the total curvature of $\C_2$:
$$
\int_{\C_1} \kappa_{\C_1}\,ds \le \int_{\C_2}\kappa_{\C_2} \,ds.
$$
\end{proposition}

\begin{proof}
This is obvious from Figure~\ref{fig:inside} and the interpolation of the
total curvature as the change in angle along the curve.  The reader wanting a
more detailed (or more highbrow) proof can construct one from the Gauss-Bonnet
formula for surfaces with boundaries having canners.  For this see
Equation~(4) in the
excellent  expository article \cite{Chern:triangles} by S.-S.
Chern.
\end{proof}

Various elementary inequalities between bounds on the length, total curvature,
and bounds on the radius will come up often enough that it is worth
recording them.

\begin{proposition}\label{prop:tc-len-bds}
If the radius of curvature of a curve satisfies $R_1\le \rho \le R_2$
for some positive constants $R_1$ and $R_2$ and if $L$ is the
length of $\C$ then
$$
R_1 \tc(\C) \le L \le R_2 \tc(\C)
$$
and
$$
\frac{L}{R_1^{1/3}} \le \left( \frac{\tc(\C)R_2}{R_1} \right)^{1/3} L^{2/3}
\le \frac{\tc(\C)R_2}{R_1^{1/3}} .
$$
\end{proposition}

\begin{proof}
The first of these follows from $\tc(\C) = \int_\C \kappa\,ds= \int_\C (1/\rho)\,ds$
and $R_1\le \rho\le R_2$.  The second follows from just using $L \le \tc(\C)R_2$
\begin{align*}
	\frac{L}{R_1^{1/3}} &= \left( \frac{L}{R_1} \right)^{1/3}L^{2/3}
	\le \left( \frac{\tc(\C)R_2}{R_1} \right)^{1/3} L^{2/3}\\
	&\le \left( \frac{\tc(\C)R_2}{R_1} \right)^{1/3} \left( \tc(\C)R_2\right)^{2/3}
	=  \frac{\tc(\C)R_2}{R_1^{1/3}}
\end{align*}
\end{proof}

Another basic tool we will use is an elementary
maximum principle.  This is well-known, but we include a short proof
for completeness.

\begin{proposition}[Maximum Principle]\label{max_prin}
Let $\C_1$ and $\C_2$ be convex curves with $\C_1$ inside the convex
hull of $\C_2$ and tangent to $\C_2$ at some point $P$ (which could
be endpoints of $\C_1$ and $\C_2$).
Then at $P$
$$
\kappa_{\C_1} \ge \kappa_{\C_2}.
$$
\end{proposition}

\begin{figure}[h]
	\begin{overpic}[height=.4in]{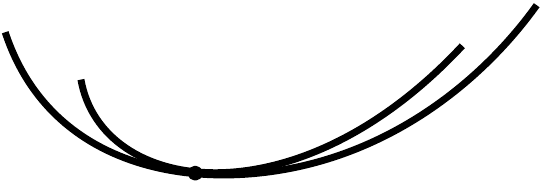}
\put(15,20){$\C_1$}
\put(0,30){$\C_2$}
\put(35,5){$P$}
\end{overpic}
\caption{$C_1$ is at least as curved as $\C_2$ at the point
$P$.}
\label{fig:above}
\end{figure}

An informal restatement is that if $\C_1$ is internally tangent to $\C_2$
at $P$, then $\C_1$ is as least as curved as $\C_2$ at $P$.

\begin{proof}
In an appropriate
coordinate system, and possibly working with smaller pieces
of the curves near $P$, we can write $\C_1$ and $\C_2$ as graphs
$y=f_1(x)$ and $y=f_2(x)$ respectively.  Then the hypothesis of the proposition
is that the function $f_1-f_2$ has a local minimum at $P$.  The first
and second derivative tests yield that if $P=(x_0,f(x_0))$,
then $f'_1(x_0)-f'_2(x_0)=0$ and $f_1''(x_0)-f_2''(x_0)\ge 0$.  Using
this equality and inequality and  the standard formula for the curvature of graphs we get
$$
\kappa_{\C_1}(P) = \frac{f_1''(x_0)}{(1+f_1'(x_0)^2)^{3/2} }
\ge  \frac{f_2''(x_0)}{(1+f_2'(x_0)^2)^{3/2} }=\kappa_{\C_2}(P).
$$
\end{proof}

Another well known fact is that two $C^2$ curves with common endpoints, tangent,
and curvature can be joined together to form a $C^2$ curve. Again, we include
a short proof.

\begin{figure}[h]
	\begin{overpic}[height=.4in]{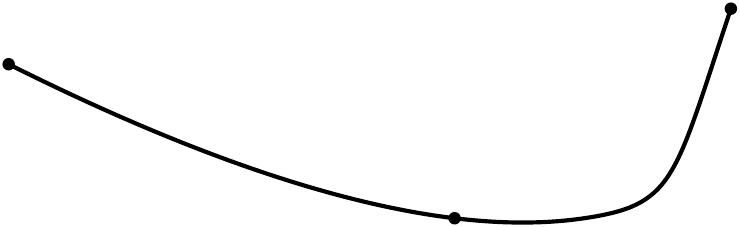}
		\put(30,10){$\C_1$}
		\put(78,7){$\C_2$}
		\put(58,3){$P$}
	\end{overpic}
	\caption{The curves  $\C_1$ and $\C_2$ have the same tangent
	and curvature at $P$. This implies the union $\C = \C_1\cup \C_2$
	is also a $C^2$ curve.}
	\label{fig:splice}
\end{figure}

\begin{lemma}[Splicing Lemma]\label{lem:splice}
Let $\C_1$ and $\C_2$ be two curves of class $C^2$ such that the
terminal point of $\C_1$ is the initial point of $\C_2$, and
that at this common point the two curves have the same tangent
and curvature as in Figure \ref{fig:splice}.  Then $\C_1\cup \C_2$ is a curve of class $C^2$.
\end{lemma}

\begin{proof}
There are coordinates so that near $P$,   $\C_1$ is the graph of a function
$y=f_1(x)$ for on the interval $[a,0]$ and $\C_2$ is the graph of $y=f_2(x)$
on $[0,b]$ for $C^2$ functions $f_1$ and $f_2$. As the terminal point of $\C_1$
is the initial point of $\C_2$ we have $f_1(0)=f_2(0)$.  That the
curves have the same tangent at this point implies $f'_1(0)=f_2'(0)$.
The equality of the curvatures at $x=0$ gives
$f_1''(0)/(1+f_1'(0)^2)^{3/2} = f_2''(0)/(1+f_2'(0)^2)^{3/2} $
which implies $f_1''(0)=f_2''(0)$.  Therefore,  the function
given by $f(x)=f_1(x)$ on $[a,0]$ and $f(x)=f_2(x)$ on $[0,b]$
is continuous with continuous first and second derivatives.
Whence $\C=\C_1\cup \C_2$ is the graph of a $C^2$ function
near the common endpoint, showing that $\C$ is $C^2$.
\end{proof}

\begin{lemma}\label{lem:lunes}
Let $\C_1$ and $\C_2$ be $C^2$ convex curves with the same endpoints
and with $\C_1$ contained in the convex hull of $\C_2$.  Assume
the total curvature of $\C_2$ satisfies
$$
\int_{\C_2}\kappa\,ds \le \pi.
$$
Then $\C_2$ is at least as curved as $\C_1$ in the sense that
$$
\max_{P\in \C_2} \kappa_{\C_2}(P) \ge \min_{Q \in \C_1} \kappa_{\C_1}(Q).
$$

\end{lemma}

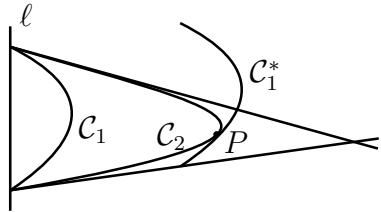
\begin{figure}[h]
\begin{overpic}[height=1in]{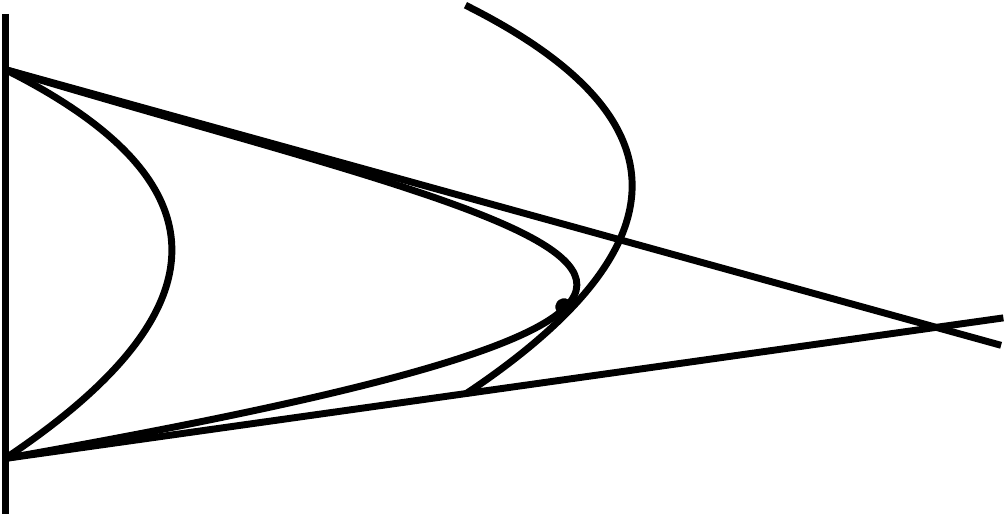}
\put(19,20){$\C_1$}
\put(65,35){$\C_1^*$}
\put(40,18){$\C_2$}
\put(58,16){$P$}
\put(3,50){$\ell$}
\end{overpic}
\caption{$\C_1$ can be translated to a position $\C_1^*$ where it
is externally tangent to $\C_2$ at the point $P$.  At $P$ the curve
$\C_2$ is at least as curved as $\C_1^*$.}
\label{fig:lunes}
\end{figure}

\begin{proof}
Let $\ell$ be the line through the common endpoints of the two curves
and consider the tangent lines to $\C_2$ at its endpoints as in Figure~\ref{fig:lunes}.
Because
the total curvature of $\C_2$ is at most $\pi$ these lines will either be
be parallel (when the total curvature is $\pi$) or will intersect
on the same side of $\ell$ as $\C_1$ and $\C_2$.  Translate $\C_1$
keeping one of its endpoints on one of the tangent lines to a position
$\C_1^*$ where it is tangent to $\C_2$ at a point $P$ (this is the farthest translated
position where $\C_2$ and $\C_1^*$ still intersect).  By the
maximum principle $\kappa_{\C_2}(P)\ge \kappa_{\C_1^*}(P)$.  As
translation preserves curvature this completes the proof.
\end{proof}

\begin{figure}[h]
\begin{overpic}[height=.5in]{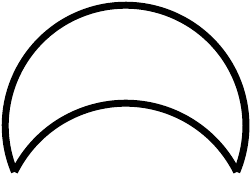}
\put(45,32){$\C_1$}
\put(45,72){$\C_2$}
\end{overpic}
\caption{The curvature of $\C_2$ is everywhere less than the curvature
of $\C_1$.}
\label{fig:circ_arcs}
\end{figure}

The hypothesis $\int_{\C_2}\kappa\,ds\le \pi$ is necessary as can be seen in the example of
the two circular arcs in Figure~\ref{fig:circ_arcs}.

\begin{theorem}\label{thm:2closed}
Let $\C_1$ and $\C_2$ be $C^2$ closed convex curves that intersect in
three or more points.  Then, they have comparable curvature in
the sense that there are points $P$ on $\C_1$ and $Q$  on $\C_2$
with $\kappa_{\C_1}(P) = \kappa_{\C_2}(Q)$.
\end{theorem}

\begin{corollary}\label{cor:close_cap_circ}
Let $\C$ be a closed convex $C^2$ curve that intersects a circle
of radius $R$ in three or more points.  Then, there is a point
on $\C$ with $\kappa=1/R$.
\end{corollary}

\begin{proof}[Proof of Theorem \ref{thm:2closed}]
If $\C_1$ and $\C_2$ intersect in infinitely many points, then
let $P$ be an accumulation point of the set of intersection points.
At $P$ the two curves will have contact of order at least $2$ (and
in fact,  infinite order if the curves are of class $C^\infty$) and
therefore,  have the same curvature at $P$.
Thus, we can assume the two curves only intersect in finitely many points.

\begin{claim}
There are points $P_j$ on $\C_j$ for $j=1,2$ such that
$\kappa_{\C_1}(P_1) \le \kappa_{\C_2}(P_2)$.
\end{claim}

Assuming the claim the theorem follows.  For the claim implies
the function $\kappa_{\C_2} - \kappa_{\C_1}$ is non-negative
at some point on the Cartesian product $\C_1\times \C_2$.  By symmetry this function
is also non-positive at some point.  As $\C_1\times \C_2$ is
connected this implies $\kappa_{\C_2} - \kappa_{\C_1}=0$ at some point,
which is equivalent to the conclusion of the theorem.

The proof of the claim splits into three cases.

\emph{Case 1: $\C_1$ is externally tangent to $\C_2$ at some point
of intersection.}  Then the claim follows directly from the
maximum principle (Proposition~\ref{max_prin}).

\emph{Case 2: $\C_1$ is internally tangent to $\C_2$ at some point of
intersection.} Let $\C_1$ be internally tangent to $\C_2$ at the point
$P$.  Let $P_-$ and $P_+$ be the points of intersection that
are on either side of $P$ (these exist as there are only finitely many
points of intersection). As the total curvature of $\C_2$ is $2\pi$, at least
on of the two arcs $\C_2\big|_{P_-}^P$ or $\C_2\big|_{P}^{P_+}$ will have
total curvature $\le \pi$.  Then Lemma \ref{lem:lunes} implies the
conclusion of the claim holds.

\emph{Case 3: At every point of intersection $\C_1$ crosses $\C_2$.}
Between each two consecutive points of intersection the arc of $\C_2$
between these points is either inside of $\C_1$, call such arcs \emph{positive},
or outside of $\C_1$, call such arcs \emph{negative}.  In the current case
each point of intersection is between a positive and negative arc of $\C_2$.
Therefore,  the total number of points of intersection is even and the number
of positive arcs of $\C_2$ is half of this number.  The number of points of
intersection is at least $3$ and therefore the $\C_2$ has
at least two positive arcs.  And again,  as the total curvature of $\C_2$
is $2\pi$ at least one of these arcs has total curvature $\le \pi$,
and again we can use Lemma \ref{lem:lunes} to see the claim holds.
\end{proof}

\begin{theorem}\label{thm:lun_circ_intersec}
 Let $\C$ be a $C^2$ convex curve with total curvature satisfying
$
\int_\C \kappa\,ds \le \pi
$
that intersects a  circle of radius $R$ in three or more points.
Then, there is a point on $\C$ with curvature $\kappa=1/R$.
\end{theorem}

\begin{proof}
We first consider the case when $\int_\C \kappa\,ds<\pi$.  Let $P_0$
be the initial point of $\C$ and $P_1$ the terminal point.

Let $\kappa_0$ be
the curvature of $\C$ at $P_0$ and $\kappa_1$ its curvature at $P_1$.
Let $\alpha = \pi -\int_\C \kappa\,ds$.  Construct a
curve $\C_1$ with total curvature $\alpha$ and with curvature $\kappa_1$
at its initial point and $\kappa_0$ at its terminal point and
with its curvature everywhere between $\kappa_0$ and $\kappa_1$.  As explicit example
of such a curve can be constructed by letting $\theta\cn [0,\alpha]\to \R$ be a function
with derivative
$$
\theta'(t) = \frac{\alpha-t}{\alpha} \kappa_1+ \frac{t}{\alpha} \kappa_0.
$$
and letting
$$
\gamma(s) =  \int_0^s ( \cos\theta(t),\sin\theta(t))\,dt.
$$
Then $\gamma$ is unit speed curve with curvature $\kappa(s) = \theta'(s)$.
By rotating and translating $\C_1$ we can move it until its initial point
is $P_1$ and $\C$ and $\C_1$ have the same tangent vector at $P_1$ as in
Figure \ref{fig:extend}.

\begin{figure}[h]
	\begin{overpic}[height=1.3in]{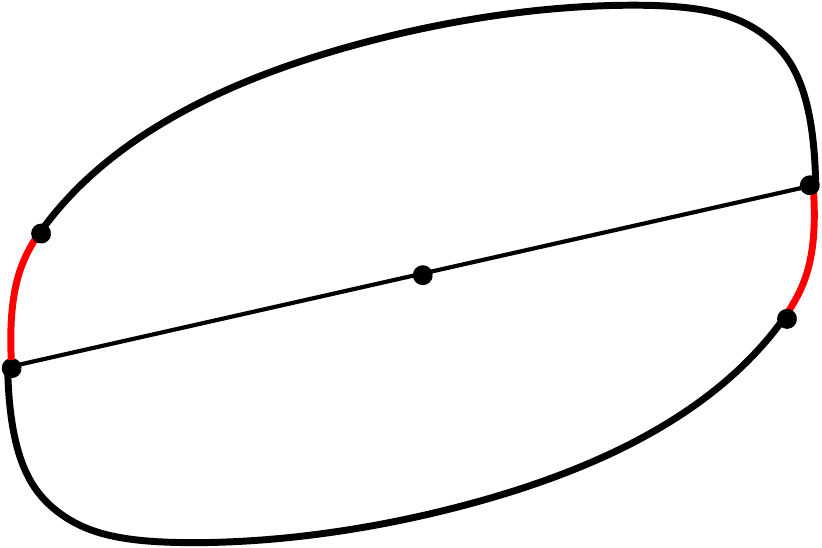}
		\put(22,2){$\C$}
		\put(65,59){$\C^*$}
		\put(3,27){$\C_1^*$}
		\put(90,34){$\C_1$}
		\put(48,24){$M$}
		\put(3,14){$P_0$}
		\put(96,21){$P_1$}
		\put(101,42){$P_2$}
		\put(7,35){$P_1^*$}
	\end{overpic}
	\caption{Extend the curve $\C$ by the curve $\C_1$ so that
	that total curvature of $\C\cup \C_1$ is $\pi$ and is
	of class $C^2$. Let $P_2$
	be the terminal point of this union.  Rotate these curves around
	the midpoint, $M$, of the segment between $P_0$ and $P_2$.  The
	resulting closed curve will be of class $C^2$ as can be seen
	by four applications of The Splicing Lemma \ref{lem:splice}.}
	\label{fig:extend}
\end{figure}
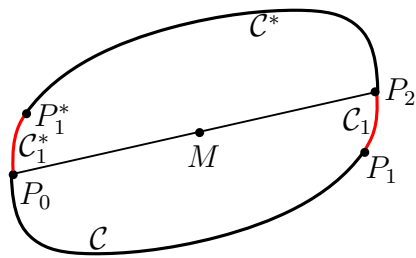

 Take the resulting curve $\C\cup \C_1$
and rotate it about the midpoint, $M$, of the segment between $P_0$ and $P_1$
and let $\C^*$ and $\C_1^*$.  Then the union $\mathcal B=\C\cup \C_1\cup \C^*\cup \C_1^*$
is a closed convex curve.  As $\C$ and $\C_1$ are $C^2$ the curve $\mathcal B$ is
$C^2$ except possibly at the points $P_0$, $P_1$, $P_2$, and $P_1^*$.
At $P_1$ the curves $\C$ and $\C_1$ have the same tangent vector and by
construction they have the same curvature at $P_1$.  Therefore $\mathcal B$
is $C^2$ in a neighborhood of $P_1$ by the Spicing Lemma \ref{lem:splice}.
A similar argument shows $\mathcal B$
is $C^2$ near the remaining points $P_0$, $P_2$, and $P_1^*$.

As $\C$ intersects some circle of radius $R$ in three or more points, the
curve $\mathcal B$ will also meet this circle in three or more points.
By Corollary \ref{cor:close_cap_circ} the curve $\mathcal B$ contains
a point $P$ where $\kappa=1/R$.  If $P$ is on $\C$ we are done.  It
$P$ is on $\C^*$, then, as $\C^*$ is just a rotation of $\C$, there
is a point of $\C$ with $\kappa=1/R$.  If $P$ is on $\C_1$,
then by the construction of $\C_1$ we have $\kappa_{\C_1}(P)=1/R$ is
between $\kappa_0$ and $\kappa_1$ and by the intermediate value theorem
there is a point  of $\C$ with curvature $1/R$.  A similar argument works in
the case when $P$ is on $\C_1^*$.  This covers all the cases and completes
the proof in the case the total curvature of $\C$ is less than~$\pi$.
\medskip

If the total curvature of $\C$ is $\pi$ and $\C$ intersects the circle of
radius $R$ in four or more points, then it will have proper sub-arc that
intersects the circle in three or more points and such that this sub-arc will
have total curvature less than $\pi$ and we are back in the case we have
just covered.  So,  assume $\C$ intersects the circle of radius $R$ in exactly
three points.  If one of the endpoints of $\C$ is not a point of intersection,
then there is again a proper sub-arc of $\C$ that contains the three points
of intersection with the circle
and this sub-arc will have total curvature less than $\pi$
and we are done.

Therefore, we can assume that $\C$ intersects the circle of radius $R$ in
exactly three points $P_0$, $P_1$, and $P_2$ and that $P_0$ and $P_2$
are endpoints of $\C$.  As the total curvature of $\C$ is $\pi$ the tangent
lines to $\C$ at $P_0$ and $P_2$ are parallel.  By a rotation we can assume
these are vertical and that $\C$ is the graph of a convex function.
The proof now splits into cases. Let $\mathcal S$ be the circle of radius
$R$ intersecting $\C$ in the points $P_0$, $P_1$, and $P_2$.

\emph{Case 1: The points $P_0$, $P_1$, and $P_2$ are all on the
closed lower half of~$\mathcal S$.}

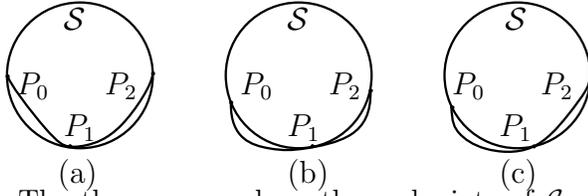
\begin{figure}[h]
	\begin{overpic}[height=.8in]{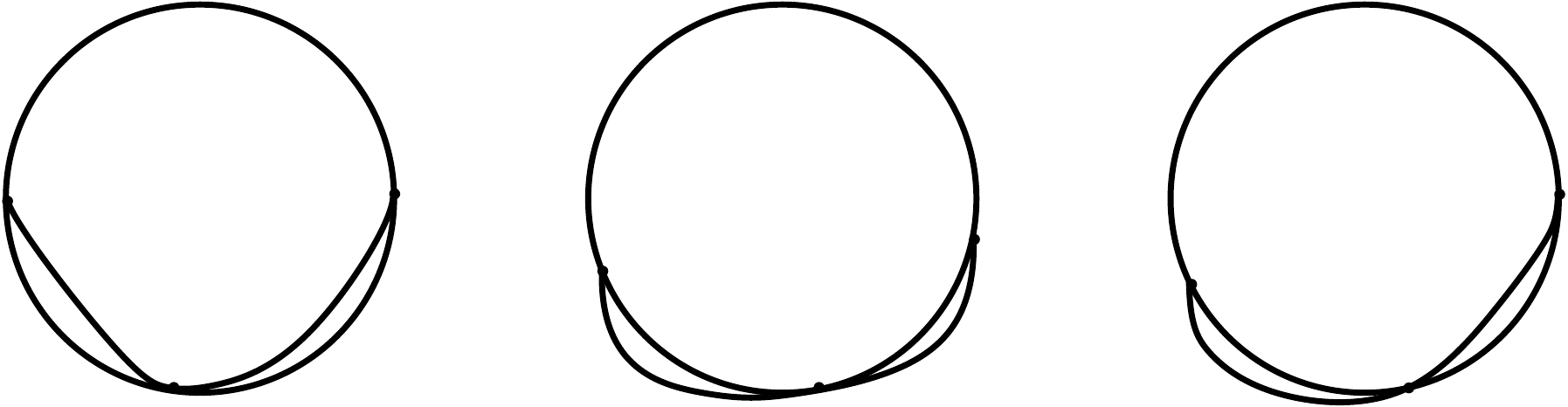}
		\put(10,21){$\mathcal S$}	
		\put(48,21){$\mathcal S$}	
		\put(85,21){$\mathcal S$}	
		\put(9,-5){(a)}	
		\put(48,-5){(b)}	
		\put(84,-5){(c)}	
		\put(2,10){$P_0$}	
		\put(40,10){$P_0$}	
		\put(76,10){$P_0$}	
		\put(10,3){$P_1$}	
		\put(48,3){$P_1$}	
		\put(84,3){$P_1$}	
		\put(17,10){$P_2$}	
		\put(55,10){$P_2$}	
		\put(91,10){$P_2$}	
	\end{overpic}
	\caption{The three cases where the endpoints of $\C$ are on the closed
	lower half of the circle $\mathcal S$.}
	\label{fig:lower}
\end{figure}

There are three sub-cases.  First, $\C$ could be internally tangent
to $\mathcal S$ at $P_1$ as in Figure \ref{fig:lower}~(a).
Then, by the Maximum Principle $\kappa_\C(P_1)\ge
1/R$.  The total curvature of the lower half circle is $\pi$ and thus
by Lemma \ref{lem:lunes} there is a point $Q$ of $\C$ between $P_0$ and $P_1$
with $\kappa_\C(Q)\le 1/R$.  Thus, there is a point on $\C$ with curvature $1/R$.

Sub-case (b) is as in Figure \ref{fig:lower}~(b) where $\C$ is externally
tangent to $\mathcal S$ at $P_1$.  By the Maximum Principle $\kappa_\C(P_1)\le 1/R$, 
and as the total curvature of $\C$ is $\pi$ Lemma \ref{lem:lunes}
gives a point between $P_0$ and $P_1$ where $\kappa_\C\ge 1/R$.  Thus,
there is a point with $\kappa_\C=1/R$.

Sub-case (c) is as in Figure \ref{fig:lower} (c) where $\C$ crosses $\mathcal S$
at $P_1$.  Then $\C$ and the lower half of the circle have total curvature
$\pi$ and therefore Lemma \ref{lem:lunes} can be applied twice, once
between $P_0$ and $P_1$ to find a point of $\C$ with $\kappa_\C\ge 1/R$,
and once between $P_1$ and $P_2$ to find a point on $\C$ with $\kappa_\C\le 1/R$.
So again, there is a point with $\kappa_\C=1/R$.

\emph{Case 2: At least one of the endpoints of $\C$ is in the open
upper half of the circle $\mathcal S$.}
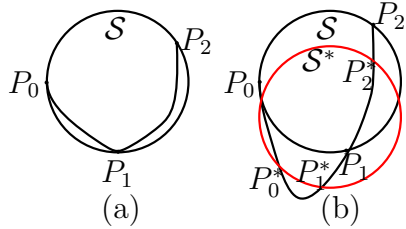
\begin{figure}[h]
	\begin{overpic}[height=1in]{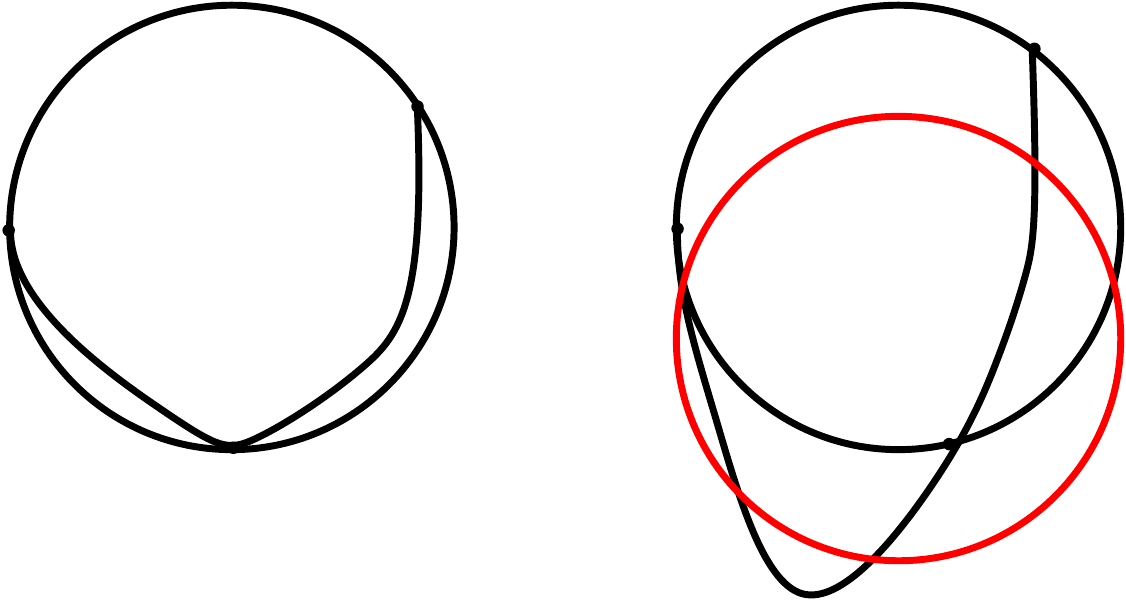}
		\put(17,45){$\mathcal S$}	
		\put(76,45){$\mathcal S$}	
		\put(72,36){$\mathcal S^*$}	
		\put(16,-5){(a)}	
		\put(77,-5){(b)}	
		\put(-10,30){$P_0$}	
		\put(50,30){$P_0$}	
		\put(16,5){$P_1$}	
		\put(82,7){$P_1$}	
		\put(38,42){$P_2$}	
		\put(92,49){$P_2$}	
		\put(83,33){$P_2^*$}	
		\put(70,5){$P_1^*$}	
		\put(57,3){$P_0^*$}	
	
	\end{overpic}
	\caption{The two cases where one endpoint of $C$ is on the open upper
	half of the circle $\mathcal S$.  }
	\label{fig:upper}
\end{figure}

Let $P_2$ be an endpoint that is in the upper half of $\mathcal S$.
As the tangent to $\C$ at $P_2$ is vertical the curve $\C$ will contain points
in the interior of $\mathcal S$.  Thus,  there are two sub-cases.

Sub-case (a) is when $\C$ is internally tangent to $\mathcal S$ at $P_1$
as in Figure \ref{fig:upper}~(a).
By the maximum principle $\kappa_\C\ge 1/R$ at $P_1$.
The total curvature of the circle is $2\pi$ and therefore at least
one of the arc $\mathcal S\big|_{P_0}^{P_1}$ or $\mathcal S\big|_{P_1}^{P_2}$
has total curvature $\le \pi$.  Lemma \ref{lem:lunes} then gives
a point of $\C$ with $\kappa_\C \le 1/R$ and there is a point with
$\kappa_\C = 1/R$.

Sub-case (b) is when $\C$ crosses $\mathcal S$ at $P_1$. As the tangent to
$\C$ at $P_2$ is vertical the part of $\C$ near $P_2$ is interior to $\mathcal S$.
Translate $\mathcal S$ downward to a position $\mathcal S^*$
so that it contains points in the region
bounded by the two curves $\mathcal S\big|_{P_0}^{P_1}$ and $\C\big|_{P_0}^{P_1}$.
As the lower half of $\mathcal S^*$ contains both points inside and outside
this region and it does not intersect $\mathcal S$, we see the lower half
of $\mathcal S^*$ intersects $\C\big|_{P_0}^{P_1}$ in at least
two points $P_0^*$ and $P_1^*$.  As $P_1$ is inside of $\mathcal S^*$
and $P_2$ is outside of $\mathcal S$ the circle $\mathcal S^*$ will
intersect $\mathcal S^*$ at some point $P_2^*$ on $\C\big|_{P_1}^{P_2}$.
Therefore $\C\big|_{P_0^*}^{P_2^*}$ intersects the circle $\mathcal S^*$
of radius $R$ in at least three points and as it is a proper sub-arc of $\C$
it has total curvature $<\pi$.  Therefore  $\C\big|_{P_0^*}^{P_2^*}$,
and thus also $\C$, has a point with curvature $=1/R$.
\end{proof}

\begin{figure}[h]
	\begin{overpic}[width=3.5in]{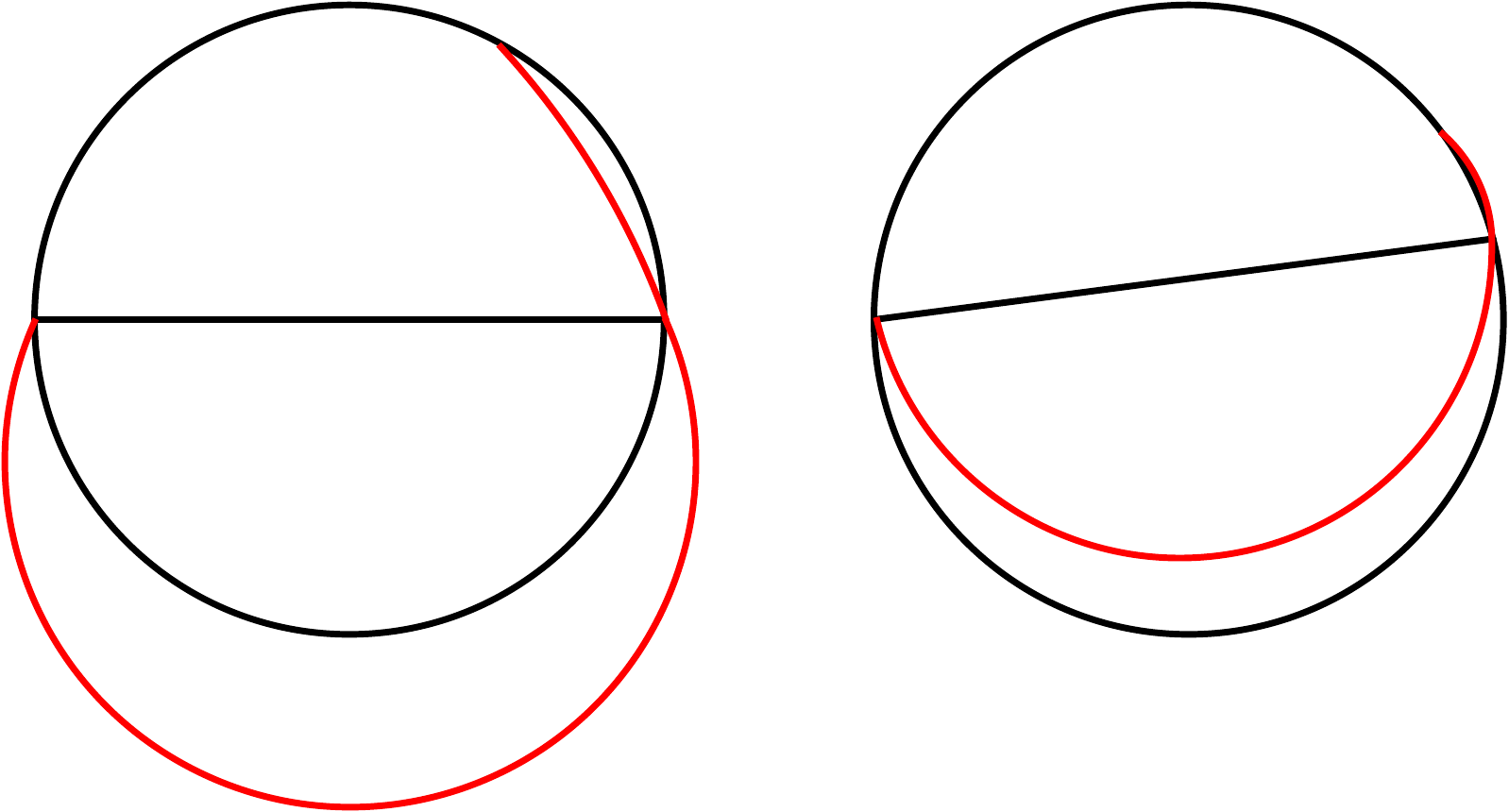}
		\put(7,7){$\C_1$}
		\put(69,19){$\C_2$}
		\put(8,45){$\mathcal S$}
		\put(64,45){$\mathcal S$}
	\end{overpic}
	\caption{The circle $\mathcal S$ has radius $R$.  In the figure
		on the left, $\C_1$ meets $\mathcal S$ in three
		points and has curvature $< 1/R$ at all points.
		In the figure on the right, $\C_2$ meets $\mathcal S$ in
		three points and has curvature $> 1/R$ at all  points.
	The total curvatures of $\C_1$ and $\C_2$ can be made arbitrarily close
to $\pi$.}
\label{fig:3pt_example}
\end{figure}
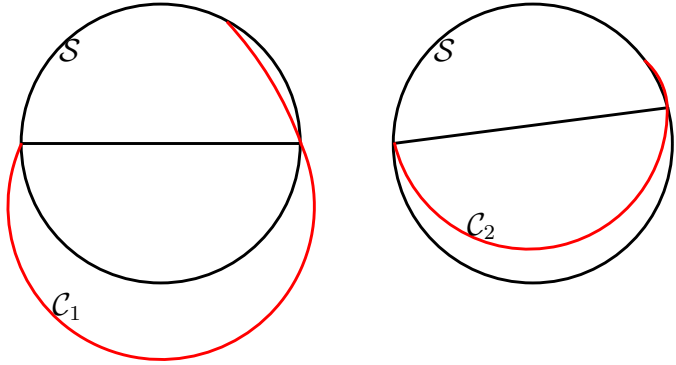

The curves in Figure~\ref{fig:3pt_example} show the hypothesis
$\int_\C \kappa\,ds\le \pi$ in
Theorem~\ref{thm:lun_circ_intersec} is best possible.

\section{Affine arclength and bounding the number of lattice points on circles and ellipses.}
\label{sec:ellipses}

\begin{theorem}\label{thm:circ<3}
	Let $\C$ be an arc of length $L$ of a circle with radius $R$
	and let $\lat$ be  a lattice.  If
	$$
	\frac{L}{R^{1/3}} \le 2 (A_\lat)^{1/3}
	$$
	then $\C$ contains at most $2$ points of the lattice $\lat$.
\end{theorem}

\begin{theorem}\label{thm:circ_open}
	Let $\C$ be an arc of length $L$ of a circle  with radius $R$
	and let $\lat$ be  a lattice. Then
	$$
	\#(\C\cap \lat) < 2 + \frac{L}{(A_\lat R)^{1/3}}.
	$$
\end{theorem}

\begin{theorem}\label{thm:circ_closed}
	Let $\mathcal S$ be a circle of radius $R$ and $\lat$ a lattice.  Then
	$$
	\#(\mathcal S\cap \lat)< \frac{\len(\mathcal S)}{(A_\lat R)^{1/3}}
	= \frac{2\pi R^{2/3}}{A_\lat^{1/3}} .
	$$
\end{theorem}

\begin{proof}[Proof of Theorem \ref{thm:circ<3}.]
If $\C \cap \lat$ has three or more points then let $P_0$, $P_1$, and
$P_2$ be distinct points  in $\C \cap \lat$.  Then
the triangle $\triangle P_0P_1P_2$ has area $\ge \frac{1}{2} A_\lat$ by
Proposition \ref{A_lat-est} and    Proposition \ref{prop:Area_pts_circle} implies
$$
\frac{A_\lat}{2} < \frac{(\|P_1-P_0\|+ \|P_2-P_1\|)^3}{16 R} < \frac{L^3}{16R}
$$
which simplifies to $L> 2(A_\lat R)^{1/3}$.  This proves the contrapositive of
the theorem.
\end{proof}

\begin{proof}[Proof of Theorem \ref{thm:circ_open}.]   Let  $P_1, P_1 ,\ldots, P_N$
	be the points of $\C \cap \lat$.  Then Proposition \ref{A_lat-est}
	implies the hypothesis of Theorem \ref{basic_est_open} holds with
	$A_0=A_\lat$ and $R_0=R$.
\end{proof}

\begin{proof}[Proof of Theorem \ref{thm:circ_closed}.]
	This proof is identical to the previous proof except that this time
	Theorem \ref{basic_est_closed} rather than Theorem \ref{basic_est_open} is
	used.
\end{proof}

\begin{remark}\label{2_is_best}
The constant $2$ in Theorem \ref{thm:circ<3} is sharp.  Let $\mathcal S$ be
a circle of radius $R$ and let $P_0$, $P_1$, and $P_2$ be three points
on $\mathcal S$ with the arclength from $P_0$ to $P_1$ and the arclength from
$P_1$ to $P_2$ being $L/2$ as in Figure \ref{fig:close}.

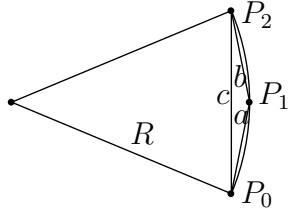
\begin{figure}[h]  
\begin{overpic}[height=1in]{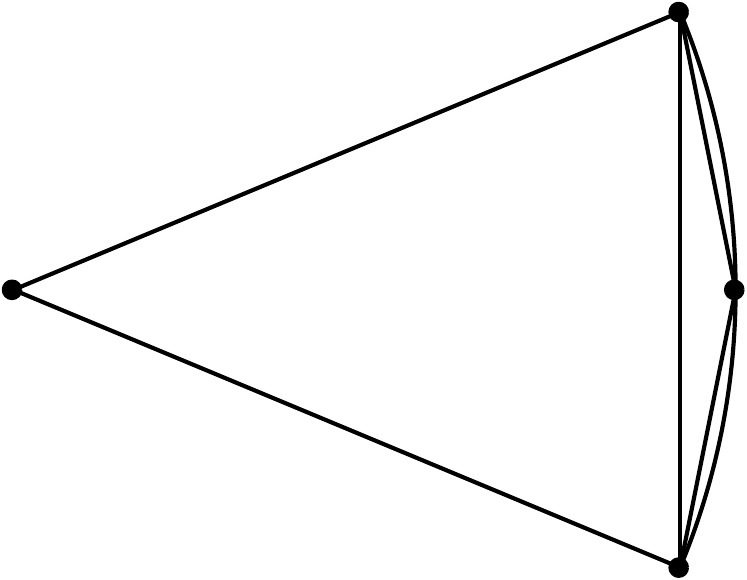}
\put(95,-2){$P_0$}
\put(95,71){$P_2$}
\put(102,38){$P_1$}
\put(50,22){$R$}
\put(85,38){$c$}
\put(92,30){$a$}
\put(92,45){$b$}
\end{overpic}
\caption{ The points $P_0$, $P_1$, and $P_2$
are on a circle of radius $R$, the arclength between $P_0$ and $P_2$ is
$L$ and $P_1$ is the midpoint of the arc between $P_0$ and $P_1$. The lengths,
$a$, $b$, and $c$ of the sides of $\triangle P_0P_1P_2$ are as shown.}
\label{fig:close}
\end{figure}

Recalling if two points $P$ and $Q$ on a circle of radius $R$
are the endpoints of an arc of length $\lambda$ on the circle,
then $\|P - Q\| = 2R\sin(\lambda/(2R))$ we find that
the side lengths of $\triangle  P_0 P_1 P_2$ are given by
$$
a(R)=b(R) = 2R \sin\left( \frac{L}{4R} \right),\qquad c(R) = 2R\sin\left( \frac{L}{2R} \right)
$$
Let $\lat_R=\lat(P_0,v_1,v_2)$ where $v_j=P_j-P_0$ for $j=1,2$.  For this lattice
$$
A_{\lat_R} = 2 \area(\triangle P_0P_1P_2) = \frac{a(R)b(R)c(R)}{2R}
$$
where the second equality follows from Proposition \ref{prop:Area_pts_circle}.
Using these in the inequality in Theorem \ref{thm:circ<3} and doing a bit of
algebra gives
$$
L \le 2 (A_{\lat_R})^{1/3} = 2 \left( \frac{a(R)b(R)c(R)}{2} \right)^{1/3}.
$$
However,
$$
\lim_{R\to \infty} a(R) = \lim_{R\to\infty} b(R)= \frac{L}{2} ,\qquad
\lim_{R\to \infty} c(R) = L.
$$
Therefore
$$
L \le \lim_{R\to \infty} 2 \left( \frac{a(R)b(R)c(R)}{2} \right)^{1/3} =
2\left( \frac{L^3}{8} \right)^{1/3} =L,
$$
showing the inequality is sharp.

This example is a bit unsatisfying as we are choosing the lattice
to depend on both $R$ and $L$.  A natural question is given
a lattice $\lat$ what is the best constant $C_\lat$ such that for any arc
of length $L$ on a circle of radius $R$ with
$$
L < C_\lat(A_\lat R)^{1/3}
$$
contains at most $2$ points of $\lat$.  Theorem \ref{thm:circ<3}
together with these examples shows
$$
\inf_{\lat} C_\lat = 2.
$$
For the lattice $\lat = \Z^2$ and restricting to circles centered
at the origin Cilleruelo \cite{Cilleruelo} and
Cilleruelo and Granville \cite{Cilleruelo-Granville}
have shown $C_{\Z^2} = 2 (2^{1/3})$.  To the best of our knowledge
$C_\lat$ is not known for any other lattice.

\end{remark}

We recall the definition of affine arclength.  Let $\C$ be
a curve with positive curvature and let $\gamma \cn [a,b]\to \C$ be a
parametrization of $\C$.  Then,  the \bi{affine arclength} of $\C$ is
given by
$$
\aff(\C) = \int_a^b (\gamma'(t)\wedge \gamma''(t))^{1/3}\,dt.
$$
If $\phi(v) = Mv+b$ is an affine map with $\det(M)>0$, then
it is straightforward to check that if $c(t) = \phi(\gamma(t))$ then
$$
c'(t)\wedge c''(t) =( M \gamma'(t)) \wedge( M \gamma''(t)) =
\det(M) \gamma'(t)\wedge\gamma''(t)
$$
and therefore,  the affine arclength transforms under affine maps with
positive determent by the rule
\begin{equation}\label{eq:aff_len_transformtion}
\aff(\phi[\C]) = \det(M)^{1/3}\aff(\C).
\end{equation}

It $\gamma\cn [a,b] \to \R^2$ is unit speed in the Euclidean sense and has
positive curvature then $\gamma'(s) = \mathbf t(s)$, and $\gamma''(s)
=\kappa(s)\mathbf n(s)$.  Thus,  $\gamma'(s)\wedge\gamma''(s)=\kappa(s)$.
This implies:

\begin{lemma}\label{Euc_to_affine}
Let $\C$ be a $C^2$ curve with positive curvature  $\kappa$.
Then,  the affine arclength of $\C$ is
\begin{equation}\label{eq:aff=kappa**1/3}
\aff(\C) = \int_a^b\kappa^{1/3}\,ds.
\end{equation}
In particular,  if  $\C$ is an arc on a circle of radius $R$, then
$$
\aff(\C) = \frac{\len(\C)}{R^{1/3}}.
$$
\end{lemma}

By an \bi{ellipse} we mean a curve $\mathcal E$ with an equation of the
form
$$
A(x-x_0)^2 + 2B(x-x_0)(y-y_0) + C(y-y_0)^2 = 1
$$
where the matrix $\left[ \begin{matrix}
A&B\\B&C
\end{matrix}\right]$ is positive definite. A fact we will use is
that if $\mathcal E$ is an ellipse, then there is affine map $\phi$ such
that $\det(\phi)=1$ and the image $\phi\big[\mathcal E\big]$ is a circle.

\begin{theorem}\label{thm:ell<3}
Let $\C$ be an arc on an ellipse $\mathcal E$ and let $\lat$ be a lattice such that
$$
\aff(\C)\le 2 (A_\lat)^{1/3}.
$$
Then,  $\C$ contains at most 2 points of the lattice $\lat$.
\end{theorem}

\begin{theorem}\label{thm:ell_open}
If $\C$ is an arc on an ellipse and $\lat$ is a lattice, then
$$
\#(\C\cap \lat) < 2 + \frac{\aff(\C)}{A_\lat^{1/3}} .
$$
\end{theorem}

\begin{theorem}\label{thm:ell_closed}
Let $\mathcal E$ be an ellipse and $\lat$ a lattice.  Then
$$
\#(\mathcal E\cap \lat) < \frac{\aff(\mathcal E)}{A_\lat^{1/3}} .
$$
\end{theorem}

\begin{proof}[Proof of Theorem \ref{thm:ell<3}]
Choose an affine map $\phi$ with $\det(\phi)=1$ and such that
the image $\mathcal S:= \phi\big[\mathcal E\big]$ is a circle
and let $R$ be the radius of this circle.
Let $\C^* = \phi\big[\C\big]$ and $\lat^*=\phi\big[\lat\big]$. 
Then by the invariance property of affine arclength under affine maps,
Lemma \ref{Euc_to_affine}, and $A_{\lat^*} = \det(M)A_{\lat} = A_{\lat}$
$$
\frac{\len(\C^*)}{R^{1/3}} = \aff(\C^*) = \det(M)^{1/3}\aff(\C)=\aff(\C)\le   (A_\lat)^{1/3}
= (A_{\lat^*})^{1/3}.
$$
Thus by Theorem \ref{thm:circ<3} $\#(\C^*\cap \lat^*)\le 2$. But $\phi$ is
a bijection so this implies $\#(\C\cap \lat)\le 2$.
\end{proof}

\begin{proof}[Proofs of Theorems \ref{thm:ell_open} and \ref{thm:ell_closed}]
	Using the notation of the proof of Theorem~\ref{thm:ell<3},
	the invariance properties of affine arclength, and
	the equalities
$\aff(\C) = \aff(\C^*) = \len(\C^*)/R^{1/3}$, $\aff(\mathcal E) = 2\pi R^{2/3}$,
$A_{\lat}=A_{\lat^*}$, $\#(\C\cap \lat) = \#(\C^*\cap \lat^*)$,
and $\#(\mathcal E\cap \lat) = \#(\mathcal S \cap \lat^*)$ hold.
Thus Theorems \ref{thm:ell_open} and \ref{thm:ell_closed} follow
directly from Theorems \ref{thm:circ_open} and \ref{thm:circ_closed}.
\end{proof}

\section{Bounding the number of lattice points on a curve by curvature and arclength.}
\label{sec:on_curve}

\begin{theorem}\label{thm:curv_bds_1}
Let $\C$ be a convex curve whose radius of curvature satisfies $\rho\ge R_1$
for some constant $R_1>0$ and whose total curvature satisfies $\int_\C\kappa\,ds\le \pi$.
Let $\lat$ be a lattice.  If
$$
\frac{\len(\C)}{(A_\lat R_1)^{1/3}}\le 2
$$
then $\C$ contains at most two points of $\lat$.
\end{theorem}

\begin{proof}
Towards a contradiction assume $\C$ contains three points $P_0$, $P_1$,
and $P_2$ of $\lat $ and that $P_1$ is between $P_0$ and $P_2$ on
$\C$.  Let $R$ be the radius of the circle through 
these points.
Because the total curvature of $\C$ is at most $\pi$
Theorem \ref{thm:lun_circ_intersec} yields a point is a point of
$\C$ with radius of curvature $R$ which implies $R\ge R_1$.
As the points $P_0$, $P_1$, and $P_2$ are in $\lat$
the lower bound $\area(\triangle P_0P_1P_2)\ge A_\lat/2$ holds
by Proposition \ref{prop:Area_pts_circle} and
$$
\frac{A_\lat}{2}\le \area(\triangle P_0P_1P_2)< \frac{(\|P_1-P_0\|+\|P_2-P_1\|)^3}{16R_1}
\le \frac{\len(\C)^3}{16R_1}
$$
which contradicts $\len(\C)/ (A_\lat R_1)^{1/3}\le 2 $.
\end{proof}

\begin{theorem}\label{thm:curv_bds_open}
Let $\C$ be an open convex curve such that the radius of convergence
of $\C$ satisfies the inequality $\rho\ge R_1$ and let $\lat $ be
a lattice.  Then
\begin{equation}\label{eq:open_curve_bd_big}
\#( \C \cap \lat)< 4 + \frac{\len(\C)}{(A_\lat R_1)^{1/3}}.
\end{equation}
If the total curvature of $\C$ satisfies $\tc(\C) \le \pi$
this can be improved to
\begin{equation}\label{eq:open_curve_bd_sm}
\#( \C \cap \lat)< 2 + \frac{\len(\C)}{(A_\lat R_1)^{1/3}}.
\end{equation}
\end{theorem}

\begin{proof}
Let $N =\#(\C\cap \lat)$ and let $P_1, P_2 ,\ldots, P_N$ be the points
of $\C \cap \lat$ listed in order along $\C$.  Let $r_j$ be
the radius of the circle through $P_j$, $P_{j+1}$, and $P_{j+2}$.
If the total curvature of $\C\big|_{P_j}^{P_{j+2}}$ is $\le \pi$
then Theorem \ref{thm:lun_circ_intersec} gives a point $Q_j$ on
this curve with $r_j =\rho(Q_j)\ge R_1$.  Also,  by Proposition \ref{A_lat-est}
$\area(\triangle P_jP_{j+1}P_{j+2})\ge A_\lat/2$.  Therefore, if
the total curvature of $\C\big|_{P_j}^{P_{j+2}}$ is $\le \pi$
for $j\in\{ 1,2 ,\ldots, N-2\}$ Theorem \ref{basic_est_open} applies
and the inequality \eqref{eq:open_curve_bd_sm} holds.  Thus will
be the case if $\tc(\C)\le \pi$.

This leaves the case where for some $k\in \{ 1,2 ,\ldots, N-2\}$
the total curvature of $\C\big|_{P_k}^{P_{k+2}}$ is greater than $\pi$.
As the total curvature of $\C$ satisfies $\tc(\C)<2\pi$
at least one of the arcs $\C\big|_{P_1}^{P_{k+1}}$ or $\C\big|_{P_{k+1}}^{P_N}$
will have total curvature $<\pi$.  We prove the case where
$\C\big|_{P_1}^{P_{k+1}}$ has total curvature $<\pi$, the other
case being similar.  Then both $\C\big|_{P_1}^{P_{k+1}}$
and $\C\big|_{P_{k+2}}^{P_{N}}$ will have total curvature $<\pi$
and by what we have just done
\begin{align*}
\#( \C\big|_{P_1}^{P_{k+1}})&< 2 + \frac{\len(\C\big|_{P_1}^{P_{k+1}} )}{(A_\lat R_1)^{1/3}}\\
\#( \C\big|_{P_{k+2}}^{P_{N}})&< 2 + \frac{\len(\C\big|_{P_{k+2}}^{P_{N}})}{(A_\lat R_1)^{1/3}}.
\end{align*}
Adding these and using $\len(\C\big|_{P_1}^{P_{k+1}} )+ \len(\C\big|_{P_{k+2}}^{P_{N}})
< \len(\C)$ shows the bound \eqref{eq:open_curve_bd_big} holds.
\end{proof}

\begin{theorem}\label{thm:curv_bds_2}
Let $\C$ be a closed convex curve whose radius of curvature satisfies
$\rho\ge R_1$ for some positive constant $R_1$ and let $\lat$ be a lattice.
Then
\begin{equation}\label{eq:close_curv-bds1}
	\#(\C\cap \lat)  <\frac{\len(\C)}{(A_\lat R_1)^{1/3}} .
\end{equation}
\end{theorem}

\begin{proof}
Let $P_1, P_2, \ldots, P_N$ be the points of $\C\cap \lat$ listed in
cyclic order around $\C$.
By Corollary \ref{cor:close_cap_circ}
the circle through $P_j$, $P_{j+1}$, and $P_{j+2}$ has radius $\rho_\C(Q)$
for some point $Q$ on $\C$ and therefore this radius is at least $R_1$.
By Proposition \ref{A_lat-est} the area of $\triangle P_jP_{j+1}P_{j+2}$ 
is at least $A_\lat/2$.  Therefore Theorem \ref{basic_est_closed}
implies \eqref{eq:close_curv-bds1}.
\end{proof}

\begin{corollary}\label{rho<=R_2}
	In Theorems \ref{thm:curv_bds_1}, \ref{thm:curv_bds_open},
	and \ref{thm:curv_bds_2} if there is also an upper bound $\rho\le R_2$
	on the radius of curvature, then the theorems still hold if
	the expression
	$$
	\frac{L}{(A_\lat R)^{1/3}}
	$$
	is replaced by either of the expressions
	$$
	\left( \frac{\tc(\C) R_2}{R_1} \right)^{1/3} L^{2/3}, \qquad
	\frac{\tc(\C)R_2}{(A_\lat R_1)^{1/3}}.
	$$
\end{corollary}

\begin{proof}
	This follows from the inequalities of Proposition \ref{prop:tc-len-bds}.
\end{proof}

\begin{remark}
The expression $\left( \frac{\tc(\C) R_2}{R_1} \right)^{1/3} L^{2/3}$
is of interest as  the coefficient $\left( \frac{\tc(\C) R_2}{R_1} \right)^{1/3}$
is invariant under dilations of the curve.  The expression
$\frac{\tc(\C)R_2}{(A_\lat R_1)^{1/3}} $ is interesting as it only depends
on the integral of curvature $\int_\C \kappa\,ds$ and the  curvature bounds
$1/R_2 \le \kappa\le 1/R_1$.
\end{remark}

\section{Bounding the number of lattice points near a curve.}
\label{sec:near_curve}

\begin{lemma} \label{lem:tri1pt}
Let $P_1,P_1,P_3,P_3'$ be points in $\R^2$. Let $\delta \geq 0$
and $\|P_3-P_3'\| \leq \delta$. Denote by $A$ the area of  $\triangle P_1P_2P_3$,
and by $A_1$ the area of  $\triangle P_1P_2P_3'$. Then,
$$
|A - A_1| \leq \frac{ \|P_1-P_2\|\delta}{2}.
$$
\end{lemma}

\begin{proof}
Let $\overleftrightarrow{P_1P_2}$ be the line through $P_1$ and $P_2$
and let $h$ be the distance of $P_3$ and $h'$ the distance of $P_3'$
from this line.  Then
$$
A = \frac{\|P_1-P_2\|h}{2} \qquad \text{and} \qquad A_1 = \frac{\| P_1-P_2\|h'}{2}.
$$
The distance between $P_3$ and $P_3'$ is at most $\delta$ and therefore $|h-h'|\le \delta$.
From this it follows
$$
|A-A_1| = \frac{\|P_1-P_2\||h-h'|}{2}\le \frac{\|P_1-P_2\|\delta}{2}.
$$
\end{proof}

\begin{lemma} \label{lem:4}
Let $\triangle P_1P_2P_3$ and $\triangle P_1'P_2'P_3'$
be triangles in the plane with areas $A$ and $A'$ respectively.
Let $\delta \geq 0$  and
assume
$$
\| P_j - P_j'\| \le \delta \quad \text{for}\quad j=0,1,2.
$$
Then
\begin{align}\label{3side_est}
|A-A'| &\le \frac{(\|P_1-P_2\| + \|P_3-P_2\| + \| P_3-P_1\|)\delta}{2}+\frac{3\delta^2}{2}\\
&\le  (\| P_2-P_1\| + \| P_3-P_2\|)\delta + \frac{ 3\delta^2}{2}
\label{2side_est}
\end{align}
\end{lemma}

\begin{proof}
Let $A=A_0$ be the area of $\triangle P_1P_2P_3$, $A_1$ the area of
$\triangle P_1P_2P_3'$, $A_2$ the area of $\triangle P_1P_2'P_3'$,
and $A_3=A'$ the area of $\triangle P_1'P_2'P_3'$.
By Lemma \ref{lem:tri1pt} and the triangle inequality
\begin{align*}
|A_0-A_1| &\le \frac{ \|P_1-P_2\|\delta}{2} \\
|A_1-A_2| & \le \frac{\|P_1-P_3'\|\delta }{2} \le \frac{(\|P_1-P_3\| +\delta)\delta }{2}\\
|A_2-A_3|&\le \frac{\|P_2'-P_3'\|\delta	}{2} \le \frac{(\|P_2-P_3\|+2\delta)\delta }{2}
\end{align*}
Therefore
\begin{align*}
|A-A'|&= |A_0-A_3|\\
&\le |A_0-A_1|+|A_1-A_2| + | A_2-A_3|\\
&\le  \frac{(\|P_1-P_2\| + \|P_3-P_2\| + \| P_3-P_1\|)\delta}{2}+\frac{3\delta^2}{2}
\end{align*}
which proves the inequality \eqref{3side_est}.  By the triangle inequality
$\| P_3-P_1\| \le \|P_2-P_1\| + \| P_3-P_2\|$ and therefore the inequality \eqref{2side_est}
follows from~\eqref{3side_est}.
\end{proof}

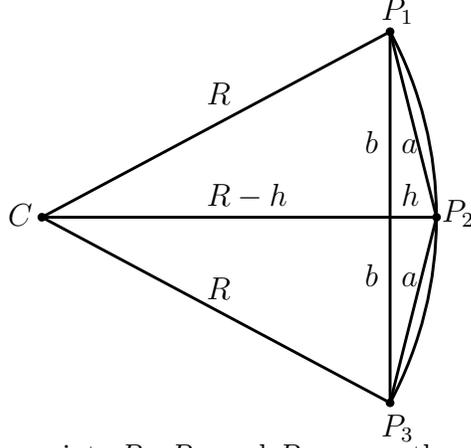
\begin{figure}[h]
	\begin{overpic}[height=2in]{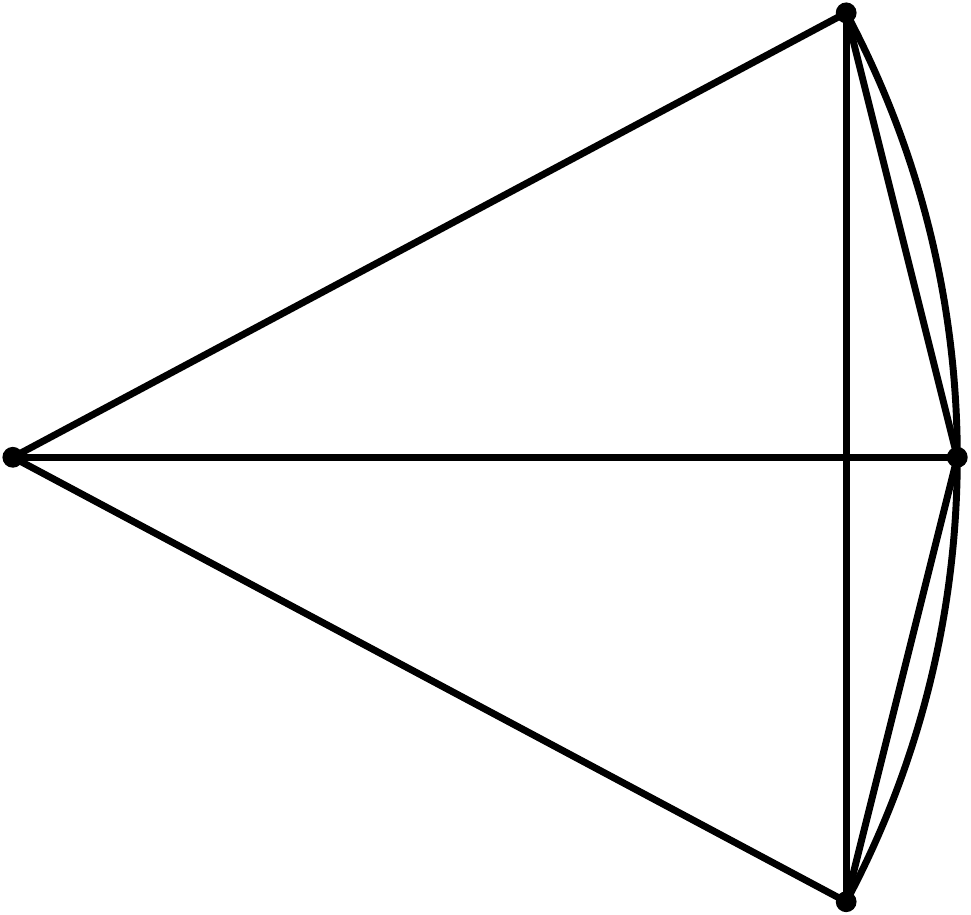}
		\put(85,96){$P_1$}
		\put(85,-7){$P_3$}
		\put(100,46){$P_2$}
		\put(42,75){$R$}
		\put(42,27){$R$}
		\put(42,50){$R-h$}
		\put(90,50){$h$}
		\put(81,63){$b$}
		\put(81,30){$b$}
		\put(90,63){$a$}
		\put(90,30){$a$}
		\put(-7,45){$C$}
	\end{overpic}
	\caption{The points $P_1$, $P_2$, and $P_3$ are on the circle of radius $R$
	centered at $C$ and $\| P_2-P_1\|=\|P_3-P_2\|=a$, $\| P_1-P_3\| =2b$,
	and $h$ is as shown.  Then $h= a^2/(2R)$.}
	\label{fig:h_formula}
\end{figure}

\begin{lemma}\label{lem:h_formula}
In the triangle $\triangle P_1P_2P_2$ as in Figure \ref{fig:h_formula}
using the side $\overline{P_1P_3}$ as a base, the height is
$$
h = \frac{a^2}{2R} .
$$

\end{lemma}

\begin{proof} 
Two applications of the Pythagorean Theorem give
$$
b^2 + (R-h)^2 = R^2, \qquad b^2+h^2= a^2
$$
Solving these for $b^2$ and setting the results equal gives
$$
R^2-(R-h)^2 = a^2- h^2
$$
and solving this for $h$ gives the desired formula.
\end{proof}

\begin{lemma}\label{lem:not_colinear}
Let $P_1,P_2,P_3$ be points on a circle of radius $R$ and let
$P_1',P_2',P_3'$  be points with $\|P_j'-P_k'\|\ge d$ 
when $j\ne k$ for some $d>0$   and $\|P_j-P_j'\|\le \delta$ for $j=1,2,3$.  Then
\begin{equation}\label{eq:delta_bd}
\delta <  \frac{d^2}{2\big(R+d+ \sqrt{(R+d)^2-d^2}\,\big)}
\end{equation}
implies the points $P_1'$, $P_2'$, and $P_3'$ are not collinear.
\end{lemma}

\begin{proof}
First note
$$
 \frac{d^2}{2\big(R+d+ \sqrt{(R+d)^2-d^2}\,\big)}
 < \frac{d^2}{2d} = \frac{d}{2}
$$
and thus the inequality \eqref{eq:delta_bd} implies
$\delta<d/2$.  Whence
$$
\|P_j-P_k\|\ge \|P_j'-P_k'\|- \|P_j-P_j'\|-\|P_k-P_k'\|> d - 2\delta>0.
$$
Therefore,  the point $P_1$, $P_2$, and $P_3$ are distinct.  

Given three points on a circle, then  at least one of the points, $P$,
is such that the other two are on opposite sides of the diameter through $P$.
For if for one of the points, call it $Q$ the other two are both on
the same side of the diameter through $Q$, or one of them is the
other end of the diameter through $Q$, then all three are on a closed
half circle.  Then let $P$ be the one of the points on this half circle
which is between the other two.  

Therefore, we can label the points $P_1$, $P_2$ and $P_3$
so that $P_1$ and $P_3$ are on opposite sides of the diameter through
$P_2$.  Farther, we can
assume the circle containing $P_1$, $P_2$, and $P_3$
goes through the origin, that $P_2$ is at the origin, the circle
is above the $x$-axis and
is tangent to the $x$-axis at $P_2$, $P_1$
is on the  open half plane defined by $x<0$, and $P_3$ is on the  
open half plane defined by $x>0$.
Let $a= d-2\delta$ and let $Q_1$ and $Q_3$  be the points  on
the circle with $\|P_2-Q_1\|=\|P_2-Q_3\|=a$ 
as in Figure~\ref{fig:not_colinear}.
Let $h$ be the distance between $P_2$ and the line
through $Q_1$ and $Q_3$.

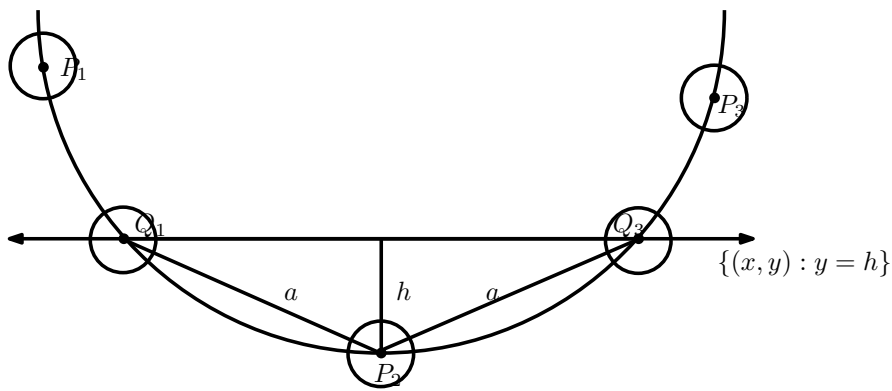
\begin{figure}[h]\footnotesize
	\begin{overpic}[height=2in]{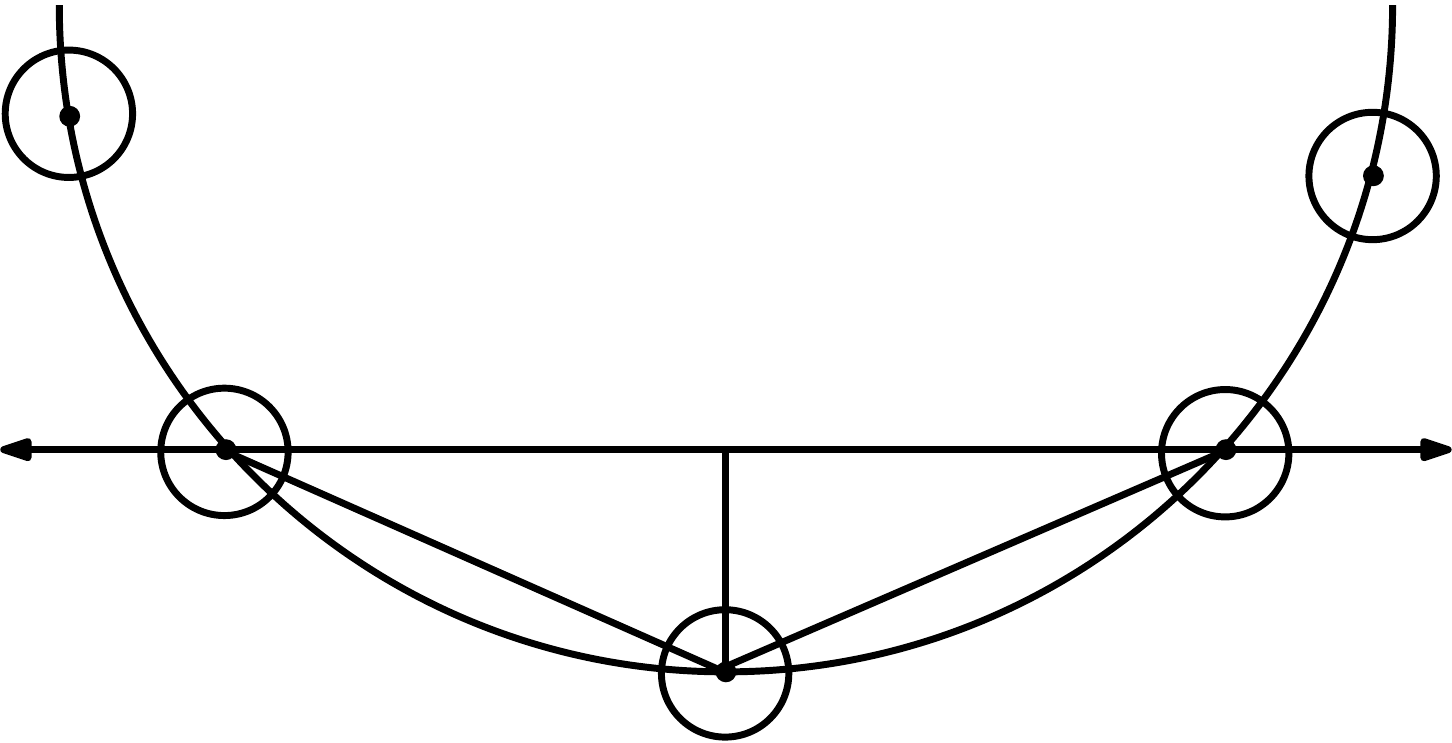}
		\put(49,1){$P_2$}
		\put(7,42){$P_1$}
		\put(95,37){$P_3$}
		\put(17,21){$Q_1$}
		\put(81,21){$Q_3$}
		\put(52,12){$h$}
		\put(64,12){$a$}
		\put(37,12){$a$}
		\put(95,16){$\{(x,y):y=h\}$}
	\end{overpic}
	\caption{The points $Q_1$ and $Q_3$ are the points on the circle
		such that $\|Q_j-P_2\|=d-2\delta=:a$.
		The circles around the points $P_j$ and $Q_j$ have radius
		$\delta$
		and therefore contain the points $P_1'$, $P_2'$,
		and $P_3'$.
		}
	\label{fig:not_colinear}
	
\end{figure}
By Lemma \ref{lem:h_formula}
$$
h = \frac{a^2}{2R} = \frac{(d-2\delta)^2}{2R}.
$$
As long at $h>2\delta$ the line with equation $y=h/2$ separates
the open disks of radius $\delta$ about $P_1$ and $P_3$ from
the open disk of radius $\delta$ about $P_2$ and therefore
the points $P_1'$ and $P_3'$ are above the line $\{y=h/2\}$
and $P_2'$ is below this line and thus these three points
are not collinear.

The inequality $h > 2\delta$ is
\begin{equation}\label{eq:2delta<h}
2\delta < \frac{(d-2\delta)^2}{2R}
\end{equation}
which is equivalent to
$$
0< 4 \delta^2 - 4(R+d)\delta +d^2.
$$
Viewing the right-hand side of this as a quadratic polynomial in $\delta$
with roots $r_1<r_2$, then
$$
r_1,r_2 = \frac{R+d\pm \sqrt{(R+d)^2-d^2}}{2} .
$$
Thus,  the inequality \eqref{eq:2delta<h} has as solution set the
union $(-\infty,r_1) \cup (r_2,\infty)$.  But $r_2>d/2$ so $(r_2,\infty)$
can be ignored.  Therefore
$$
\delta< r_1 = \frac{R+d- \sqrt{(R+d)^2-d^2}}{2}=
	\frac{d^2}{2\big(R+d+ \sqrt{(R+d)^2-d^2}\,\big)}
$$
implies the points $P_1'$, $P_2'$, and $P_3'$ are not collinear.
\end{proof}

\begin{lemma}\label{Lem:3pts_close_C}
Let $\C$ be either a closed convex curve, or a convex arc with
total curvature $\le \pi$ and assume the radius of curvature of
$\C$ satisfies
$$
 R_1 \le \rho \le R_2.
$$
Let $\lat$ be a lattice and let $P_1'$, $P_2'$, and $P_3'$
be distinct points in $\lat$ with
$$
\dist(P_j',\C)<\delta
$$
where
$$
\delta <   \frac{d_\lat^2}{2\big(R_2+d_\lat+ \sqrt{(R_2+d_\lat)^2-d_\lat^2}\,\big)}  
$$ 
and let $P_j$ be a point of $\C$ with $\| P_j-P_j'\|<\delta$ for $j=1,2,3$.
Then,
$$
\area(\triangle P_1P_2P_3)\ge \frac{A_\lat}{2} - ( \|P_2-P_1\| + \|P_3-P_2\|)\delta
-\frac32 \delta^2.
$$
\end{lemma}

\begin{proof}Let $R$ be the radius of the circle through $P_1$, $P_2$, and
	$P_3$.  Then by Corollary \ref{cor:close_cap_circ}
	or Theorem \ref{thm:lun_circ_intersec} there is a point
	on $\C$ where $\rho=R$ and thus $R_1\le R \le R_2$.
	Then by Lemma \ref{lem:not_colinear} the points $P_1'$, $P_2'$,
	and $P_3'$ are not collinear and therefore by
	Proposition \ref{A_lat-est} $\area(\triangle P_1'P_2'P_3')\ge A_\lat/2$.
	 The lower bound on the $\area(\triangle P_1P_2P_3)$
	follows from Proposition \ref{lem:tri1pt}.
\end{proof}

\begin{theorem}\label{thm:<2_near}
Let $\C$ be a convex arc with total curvature $\le \pi$ and whose
radius of curvature satisfies $R_1\le \rho\le R_2$.  Let $\lat$
be a lattice and let $\delta>0$ satisfy
$$
\delta<    \frac{d_\lat^2}{2\big(R_2+d_\lat+ \sqrt{(R_2+d_\lat)^2-d_\lat^2}\,\big)}  .
$$ 
Let $L=\len(\C)$.  Then
\begin{equation}\label{eq:<=2}
  \frac{L^3}{8R_1} + 2L\delta + 3 \delta^2\le A_\lat
\end{equation}
implies
\begin{equation}
\# \{ Q\in \lat : \dist(Q,\C) < \delta \} \le 2.
\end{equation}
\end{theorem}

\begin{proof}
We prove the contrapositive: If $\# \{ Q\in \lat : \dist(Q,\C) < \delta \} \ge 3$,
then the inequality \eqref{eq:<=2} is violated.
Assume three are three or more points
of $\lat$ at a distance less than $\delta$ form $\C$.  Then there are
$P_1$, $P_2$, $P_3$, $P_1'$, $P_2'$, and $P_3'$ that satisfy the
hypothesis of Lemma \ref{Lem:3pts_close_C}. Let $R$ be the radius of
the circle through $P_1$, $P_2$, and $P_3$.  By
Theorem \ref{thm:lun_circ_intersec} and the given bounds on
$\rho$ we have $R_1\le R\le R_2$. By  Lemma \ref{Lem:3pts_close_C}
and Proposition \ref{prop:Area_pts_circle}
\begin{align*}
\frac{A_\lat}{2} -L\delta -\frac32 \delta^2
&\le \frac{A_\lat}{2} -( \|P_2-P_1\|+\|P_3-P_2\|)\delta -\frac32\delta^2\\
&\le \area(\triangle P_1P_2P_3)\\
&< \frac{(\|P_2-P_1\|+\|P_3-P_2\|)^3}{16R} \\
&\le \frac{L^3}{16R_1}
\end{align*}
which contradicts \eqref{eq:<=2}.
\end{proof}

\begin{theorem}\label{thm:near_open}
Let $\C$ be a convex arc with total curvature at most $\pi$ with radius
of curvature bounded by $R_1\le \rho\le R_2$.  Let
$\lat$ be a lattice and $\delta>0$ with
$$
\delta <  \frac{d_\lat^2}{2( R_2 + d_\lat + \sqrt{(R_2+d_\lat)^2 - d_\lat^2})}
$$
and 
\begin{equation}\label{ineq:A_1>0}
\frac{A_\lat}{2} - L\delta -\frac32\delta^2>0.
\end{equation}
Then,
$$
\#\{Q \in \lat: \dist(\C,Q)< \delta\}
< 2 + \frac{L}{\big(R_1( A_\lat -2 L\delta - 3 \delta^2)\big)^{1/3} }
$$
where $L = \len(\C)$.
\end{theorem}

\begin{proof}
Let $N = \#\{Q \in \lat: \dist(\C,Q)< \delta\}$ and
$$
\{Q \in \lat: \dist(\C,Q)< \delta\}=\{ P_1',P_2' ,\ldots, P_N'\}
$$
and  let $P_j$ be a point of $\C$ with $\dist(P_j',\C) = \| P_j-P_j'\|$.
By Lemma \ref{Lem:3pts_close_C}
$$
\area(\triangle P_jP_{j+1}P_{j+2})\ge \frac{A_\lat}{2} -
( \|P_{j+1}-P_{j}\| + \|P_{j+2}-P_{j+1}\|)\delta
-\frac32 \delta^2.
$$
and by Theorem \ref{thm:lun_circ_intersec} the circle through $P_j$, $P_{j+1}$
and $P_{j+2}$ has radius of curvature $\rho(P)$ for some point on $\C$
and therefore its radius is $\ge R_1$.  Thus,  the result follows from
Theorem \ref{basic_est_open}.
\end{proof}

\begin{theorem}\label{thm:near_closed}
Let $\C$ be a closed convex curve  with radius
of curvature bounded by $R_1\le \rho\le R_2$.  Let
$\lat$ be a lattice and $\delta>0$ with
$$
\delta <  \frac{d_\lat^2}{2( R_2 + d_\lat + \sqrt{(R_2+d_\lat)^2 - d_\lat^2})}
$$ 
and
\begin{equation}\label{ineq:A'>0}
\frac{A_\lat}{2} - L\delta -\frac32 \delta^2>0.
\end{equation}
Then,
$$
\#\{Q \in \lat: \dist(\C,Q)< \delta\}
<  \frac{L}{\big(R_1( A_\lat -2 L\delta - 3 \delta^2)\big)^{1/3} }
$$
where $L = \len(\C)$.
\end{theorem}

\begin{proof}
	Other than using Theorem \ref{basic_est_closed}
	rather than Theorem~\ref{basic_est_open} this is
	exactly the same as Theorem \ref{thm:near_open}.
\end{proof}





We record what is the specialization of these results to circles.

\begin{corollary}\label{cor:near_circles}
Assume $\C$ is an arc of a circle of radius $R$, $\lat$ is a lattice,
and $\delta>0$ satisfies
$$
\delta<   \frac{d_\lat^2}{2(R + d_\lat
	+ \sqrt{(R+d_\lat)^2 - d_\lat^2})}.
$$ 
Let $L = \len(\C)$.  Then,
\begin{enumerate}
\item[(a)]  If $L \le \pi R$ (which is equivalent to having total curvature $\le \pi$)
	and
	$$
	\frac{L^3}{8R} + 2L\delta +  3\delta^2\le A_\lat
	$$ 
	then
	$$
	\#\{ Q \in \lat: \dist(Q,\C) < \delta\} \le 2.
	$$
\item[(b)]  If $L \le \pi R$ and
	$$
	\frac{A_\lat}{2} - 2L\delta - 3\delta^2>0,
	$$
	then
	$$
	\#\{ Q \in \lat : \dist(\C,Q)<\delta \} <
	2 +  \frac{L}{(R(A_\lat - 2L\delta - 3\delta^2))^\frac13} .
	$$
\item[(c)]  If $\C$ is the entire circle (i.e.\ $L=2\pi R$) and
	$$
	\frac{A_\lat}{2} -4\pi R \delta - 3\delta^2>0,
	$$
	then
	$$
	\#\{ Q \in \lat : \dist(\C,Q)<\delta \} < \frac{2\pi R^\frac23}{(A_\lat - 4\pi R-3\delta^2)^\frac13} .
	$$
\end{enumerate}
\end{corollary}

\section{Examples}
\label{sec:examples}

The following shows that Theorem \ref{thm:curv_bds_open} and Corollary
\ref{rho<=R_2} are close to being sharp.

\begin{theorem}\label{thm:examples}
Let $\lat$ be a lattice and $n\ge 2$ an integer.  Then
there is a convex curve $\C$ of length $L$ that contains exactly $n$ points of $\lat$,
and has lower and upper bounds
$$
R_1 = \min_{P\in \C} \rho(P), \qquad R_2=\max_{P\in \C} \rho(P)
$$
for the radius of curvature of $\C$, so that the inequalities
$$
\frac{L}{(A_\lat R_1)^{1/2} } \le \left( \frac{\tc(\C)R_2}{(A_\lat R_1)}\right)^{1/3} L^{2/3}
\le \frac{\tc(\C) R_2}{(A_\lat R_1)^{1/3}} < n+2.
$$
hold.
\end{theorem}

 \begin{proof} In light of Proposition \ref{prop:tc-len-bds}
we only need to find an example with
\begin{equation}\label{ineq:n+2}
\frac{\tc(\C) R_2}{(A_\lat R_1)^{1/3}} < n+2
\end{equation}

Let $\lat=\lat(v_0,v_1,v_2)$ where we can assume $v_1\wedge v_2>0$,
by possibly replacing $v_2$ by $-v_2$.  Then $v_1\wedge v_2=A_\lat$.
Let $a>0$ and $b = a+(n-1)$.  Define a curve $\C_a$ parametrically
$c\cn [a,b] \to \R^2$
by
$$
c_a(t) = v_0 + tv_1 + \frac{t(t+1)}{2} v_2.
$$
Each of the points $c_a(k)$ with $k = a, a+1 ,\ldots, a+(n-1)$ is
a point of $\lat$ and if $c_a(t) = v_0+ tv_1 + (t(t+1)/2)v_2 v_0=v_0 +  kv_1 + mv_2$
is a point of $\lat$ on $\C_a$, then the linear independence of $v_1$ and $v_2$
implies $k=t$ and $m = t(t+1)/2$, so that $c_a(t) = c_a(k)$.  Thus, there
are exactly $n$ points of $\lat$ on $\C_a$. We will show if $a$
is sufficiently large  this curve has the desired properties.
The derivatives of $c_a$ are
$$
c'_a(t) = v_1 + (t+1/2)v_2, \qquad c''_a(t) = v_2
$$
Then
$$
\lim_{a\to\infty} \frac{\|c_a'(t)\|}{a} = \lim_{a\to \infty}
\left\| \frac1a v_1 + \frac{t+1/2}{a} v_2\right\| = \|v_2\|
$$
and this limit holds uniformly in $t\in [a,b]$.
This gives the asymptotic formula
$$
\| c_a'(t)\| \sim a\|v_2\|
$$
and this holds uniformly for $t\in [a,b]$.
Using a standard formula for curvature
\begin{align*}
\rho&= \frac1{\kappa} = \frac{\|c_a'(t)\|^3}{c_a'(t)\wedge c_a''(t)} 
\sim \frac{a^3\|v_2\|^3}{A_\lat}
\end{align*}  
and this holds uniformly in $t\in [a,b]$. As this formula is independent
of $t$ we see that if
$R_1(a)$ and $R_2(a)$ are the minimum and maximum radius of curvature
on $\C_a$ then
$$
R_1(a) \sim R_2(a) \sim \frac{a^3\|v_2\|^3}{A_\lat} .
$$
Asymptotically the total curvature of $\C_a$ is
\begin{align*}
	\tau(a) &= \int_{\C_a} \frac{ds}{\rho}
= \int_a^b \frac{\|c_a'(t)\|\,dt}{ \left(\dfrac{\|c_a'(t)\|^3}{A_\lat} \right)} 
= A_\lat \int_a^b \frac{dt}{\|c_a'(t)\|^2} \\
&\sim A_\lat \int_a^b \frac{dt}{a^2\| v_2\|^2} 
= \frac{(n+1)A_\lat}{a^2\|v_2\|^2}
\end{align*} 
where we have used $b-a=n+1$.  
Putting these formulas together gives
$$
\frac{\tau(a) R_2(a)}{(A_\lat R_1(a))^{1/3}}\sim
\frac{\left(\dfrac{(n+1)A_\lat}{a^2\|v_2\|^2}  \right)\left( \dfrac{a^3\|v_2\|^3}{A_\lat} \right)}{\left(A_\lat \dfrac{a^3\|v_2\|^3}{A_\lat} \right)^{1/3}} = n+1.
$$
Thus for sufficiently large $a$ the inequality \eqref{ineq:n+2} holds.
\end{proof}



\end{document}